\RequirePackage{fix-cm}
\documentclass[smallextended]{svjour3}       
\smartqed  
\usepackage{amssymb}
\usepackage{amsmath}
\usepackage[affil-it]{authblk}
\usepackage{cite}
\usepackage{booktabs}
\usepackage[usenames]{color}
\usepackage[hmargin={25mm,25mm},vmargin={30mm,35mm}]{geometry}
\usepackage[utf8]{inputenc}
\usepackage[colorlinks,allcolors={blue},linkcolor={black}]{hyperref}
\usepackage{tabularx}
\usepackage{xcolor}
\usepackage{multirow}
\usepackage{colortbl}
\usepackage{array}
\usepackage{algorithm2e}
\usepackage[noend]{algorithmic}
\usepackage{paralist}
\usepackage{stackrel}
\usepackage{cancel}
\usepackage{graphicx}
\usepackage{rotating}
\usepackage{pgfplots}
\pgfplotsset{compat=1.16}
\usepackage{tikz}
\usetikzlibrary{spy}
\usetikzlibrary{arrows.meta,arrows}

\usepackage{newtxtext}
\usepackage{newtxmath}
\newcommand{\Real}{\mathbb{R}}

\newcommand{\eg}{\textit{e.g., }}
\newcommand{\ea}{\textit{et al.}}
\newcommand{\ie}{\textit{i.e., }}

\newcommand{\aT}{T}

\newcommand{\Mh}{\mathcal{M}_h}
\newcommand{\Th}{\mathcal{T}_h}
\newcommand{\Fh}{\mathcal{F}_h}

\newcommand{\Poly}[1]{\mathcal{P}^{#1}}
\newcommand{\Polyd}[2]{\mathbb{P}_{#1}^{#2}}
\DeclareMathOperator{\CARD}{card}
\newcommand{\card}[1]{\CARD(#1)}

\newcommand{\normal}{\boldsymbol{n}}

\newcommand{\FT}{\mathcal{F}_T}
\newcommand{\TF}{\mathcal{T}_F}

\newcommand{\I}{\boldsymbol{I}}

\newcommand{\xVec}{\boldsymbol{x}}
\newcommand{\uVec}{\boldsymbol{u}}

\newcommand{\wVec}{\boldsymbol{w}}

\newcommand{\vVec}{\boldsymbol{v}}

\newcommand{\fVec}{\boldsymbol{f}}

\newcommand{\iVec}{\boldsymbol{i}}
\newcommand{\jVec}{\boldsymbol{j}}
\newcommand{\eulern}{\mathrm{e}}

\newcommand{\Reynolds}{\mathrm{Re}}
\newcommand{\rbrackets}[1]{\left(#1\right)}

\newcommand{\tris}{tris}
\newcommand{\quads}{quads}

\everymath{\displaystyle}

\newcommand{\dir}{{\rm D}}
\newcommand{\neu}{{\rm N}}
\newcommand{\internal}{{\rm i}}

\newcommand{\momentum}{{\rm mnt}}
\newcommand{\mass}{{\rm cnt}}

\newtheorem{scheme}{Scheme}

\newcommand{\hhohll}{{\tt HHO-HLL}\ }
\newcommand{\hhohdiv}{{\tt HHO-Hdiv}\ }
\newcommand{\hhodivfree}{{\tt HHO-DivFree}\ }

\usepackage[normalem]{ulem}
\newcounter{corr}
\definecolor{violet}{rgb}{0.580,0.,0.827}
\newcommand{\corr}[3]{\typeout{Warning : a correction remains in page \thepage}
  \stepcounter{corr}        
	      {\color{blue}\ifmmode\text{\,\sout{\ensuremath{#1}}\,}\else\sout{#1}\fi}
              {\color{red}#2}
              {\color{violet} #3}
}

\begin{document}

\title{HHO methods for the incompressible Navier-Stokes and the incompressible Euler equations}


\author{Lorenzo Botti         \and
        {Francesco Carlo} Massa$^{\dagger}$ 
}


\institute{Lorenzo Botti \at
              Universit\`a degli Studi di Bergamo \\
              Dipartimento di Ingegneria e Scienze Applicate \\
              Tel.: +39 035 2052150\\
              Fax.: +39 035 2052310\\
              \email{lorenzo.botti@unibg.it}           
           \and
           $^{\dagger}$ \textit{corresponding author}\\
           {Francesco Carlo} Massa \at
              Universit\`a degli Studi di Bergamo \\
              Dipartimento di Ingegneria e Scienze Applicate \\
              Tel.: +39 035 2052002\\
              \email{francescocarlo.massa@unibg.it}           
}

\date{Received: date / Accepted: date}

\maketitle

\begin{abstract}
We propose two Hybrid High-Order (HHO) methods for the incompressible Navier-Stokes equations and investigate their robustness
with respect to the Reynolds number.
While both methods rely on a HHO formulation of the viscous term, the pressure-velocity coupling is fundamentally different, 
up to the point that the two approaches can be considered antithetical. 
The first method is kinetic energy preserving, meaning that the skew-symmetric discretization of the convective term
is guaranteed not to alter the kinetic energy balance. 
The approximated velocity fields exactly satisfy the divergence free constraint and 
continuity of the normal component of the velocity is weakly enforced on the mesh skeleton, leading to H-div conformity.
The second scheme relies on Godunov fluxes for pressure-velocity coupling: 
a Harten, Lax and van Leer (HLL) approximated Riemann Solver designed for cell centered formulations
is adapted to hybrid face centered formulations. The resulting numerical scheme is robust up to the inviscid limit, meaning 
that it can be applied for seeking approximate solutions of the incompressible Euler equations.
The schemes are numerically validated performing steady and unsteady two dimensional test cases and evaluating the convergence rates on $h$-refined mesh sequences. 
In addition to standard benchmark flow problems, specifically conceived test cases are conducted for studying the error behaviour when approaching the inviscid limit.

\end{abstract}

\section{Introduction}
\label{intro}
In this work we propose and numerically validate two Hybrid High-Order (HHO) methods for the Incompressible Navier-Stokes (INS) equations, governing the flow of incompressible fluids.
For the sake of simplicity, we focus on a Newtonian fluid with uniform density.
Given a polygonal or polyhedral domain $\Omega\subset\Real^d$, $d\in\{2,3\}$, with boundary $\partial \Omega$, the initial velocity field $\uVec_0: \Omega \rightarrow \Real^d$ and a finite time $t_F$ 
the incompressible Navier-Stokes problem consists in finding the velocity field $\uVec: \Omega \times (0,t_F) \rightarrow \Real^d$, and the pressure field $p: \Omega  \times (0,t_F) \rightarrow \Real$,
such that
\begin{subequations}
  \label{stokesProb}
  \begin{alignat}{2}\label{stokesProb:momentum}
    \frac{\partial \uVec}{\partial t} + \nabla \cdot \left[{(\uVec \otimes \uVec) + p \I -\nu \nabla \uVec }\right] &= \boldsymbol{f} &\qquad& \text{in $\Omega \times (0,t_F)$}, \\
    \label{stokesProb:mass}
          {\nabla \cdot \uVec} &= 0 &\qquad& \text{in $\Omega\times (0,t_F)$}, \\ 
    \uVec &= \boldsymbol{g}_\dir &\qquad& \text{on $\partial\Omega_\dir\times (0,t_F)$}. \\
    \normal \cdot \left[ p \I - \nu \nabla \uVec \right] &= \boldsymbol{g}_\neu &\qquad& \text{on $\partial\Omega_\neu\times (0,t_F)$},
  \end{alignat}
\end{subequations}
where $\normal$ denotes the unit vector normal to $\partial\Omega$ pointing out of $\Omega$, $\nu$ is the (constant) viscosity, 
$\boldsymbol{g}_\dir$ and $\boldsymbol{g}_\neu$ denote, respectively, the prescribed velocity on the Dirichlet boundary $\partial\Omega_\dir\subset\partial\Omega$ and the prescribed traction on the Neumann boundary $\partial\Omega_\neu\coloneq\partial\Omega\setminus\partial\Omega_\dir$,
while $\fVec:\Omega\to\Real^d$ is a given body force.
It is assumed in what follows that both $\partial\Omega_\dir$ and $\partial\Omega_\neu$ have non-zero $(d-1)$-dimensional Hausdorff measure (otherwise, additional closure conditions are needed).

HHO methods are gaining momentum in the field of continuum mechanics having been successfully 
applied to nonlinear elasticity problems \cite{Botti.Di-Pietro.ea:17,Abbas2018}, diffusion dominated incompressible flow problems \cite{Aghili.Boyaval.ea:15,Di-Pietro.Ern.ea:16}, 
porous-media flows \cite{Botti.Di-Pietro.ea:18} and poro-elasticity problems \cite{BOTTILPoro2021,BottiMultiPoro2021}. 
Similarly to Hybridizable Discontinuous Galerkin (HDG) methods, HHO formulations rely on Degrees Of Freedom (DOFs) associated to polynomial functions 
defined over mesh elements and mesh faces, the so called elemental and skeletal DOFs. 
Since solely skeletal DOFs are globally coupled, in particular the face based stencil 
consist of all mesh faces belonging to the boundary of the two elements sharing the face, 
elemental degrees of freedom can be eliminated by computing the Schur complement.
This procedure is known as \emph{static condensation}. 
Thanks to static condensations, HHO Jacobian matrices are sparse block matrices whose block size is driven 
by the dimension of polynomial spaces in $d{-}1$ variables when considering a $d$-dimensional flow problem. 
Accordingly, when accuracy is improved by means of higher-degree $p$-type expansions, 
the number of Jacobian matrix non-zero entries grows slower than Discontinuous Galerkin (DG) methods, 
thereby reducing the computational cost associated to matrix-vector products as well as matrices memory footprint.
From the matrix assembly viewpoint HHO methods are generally more expensive than DG methods, also due to the computational cost of static condensation,
nevertheless matrix assembly and static condensation are intrinsically parallel tasks that are expected to show optimal scalability on multicore and manycore architectures.
Besides computational efficiency considerations, HHO formulations have demonstrated to be robust with respect to mesh distortion and grading \cite{BottiHHOMG2021,Badia2021hhoCond}.
These features are of crucial importance in the context of CFD applications, where boundary layers are commonly employed to improve the resolution near wall. 

Since the pioneering works \cite{Cockburn.Shu:91,Cockburn.Shu:89,Cockburn.Lin.ea:89,Cockburn.Hou.ea:90,Cockburn.Shu:98} dating back to the late 1980s, 
DG methods have gained increased popularity in computational fluid mechanics, boosted by the 1997 landmark papers \cite{Bassi.Rebay:97,Bassi.Rebay.ea:97} on the viscous terms treatment.
In the context of incompressible flow problems, stabilized pressure-velocity was proposed by Cockburn \ea \cite{Cockburn.Kanschat.ea:02} while 
Bassi \ea \cite{Bassi.Crivellini.ea:06} exploited the artificial compressibility idea of Chorin \cite{Chorin:1967} to define suitable Godunov fluxes.
Shahbazi \ea \cite{Shahbazi:2007} proposed a semi-explicit time integration with nonlinear terms handled by means of local Lax--Friedrichs fluxes.
The extension of DG methods to general polyhedral meshes was theoretically conceived in \cite{Di-Pietro.Ern:10} and \cite{Di-Pietro.Ern:12},
leading to adaptive mesh coarsening by agglomeration \cite{Bassi.Botti.ea:12} and high-order accurate geometry representation with arbitrarily coarse meshes \cite{Bassi.Botti.ea:12*1,Bassi.Botti.ea:14}. 
$hp$-Versions and handling of meshes with small faces have been considered in \cite{Antonietti.Giani.ea:13,Antonietti.Cangiani.ea:16}; see also the recent monograph \cite{Cangiani.Dong.ea:17}.
More recently, Tavelli and Dumbser \cite{Tavelli.Dumbser:16,Tavelli.Dumbser:17} and Dumbser \ea \cite{Dumbser.ea:2018} proposed to use staggered meshes, while
Manzanero \ea \cite{Manzanero.Rubio.ea:20} devised an entropy-stable nodal dG spectral element method.

Since the seminal work of Nguyen \ea \cite{NGUYEN2011}, later analyzed by Cesmelioglu \ea \cite{Cesmelioglu.Cockburn.ea:17}, 
several HDG discretizations of the incompressible Navier-Stokes equations have been proposed. 
We mention in particular the \emph{energy stable} and momentum conserving formulation devised by Labeur \ea \cite{LabeurWellsNS2012} and the superconvergent method by Qiu \ea \cite{QiuShi2016}.
More recently, Rhebergen \ea \cite{Rhebergen.Wells:18} devised a pointwise divergence free formulation, later analysed by Kirk \ea \cite{Kirk2019}.
HHO discretizations of the Navier-Stokes equations have been originally considered in \cite{Di-Pietro.Krell:17,Di-Pietro.Krell:18}. 
A kinetic energy preserving formulation was proposed by Botti \ea \cite{Botti.Di-Pietro.ea:19*1}, see also \cite[Chapters 8 and 9]{Di-Pietro.Droniou:20} for further details.
Robustness with respect to large irrotational body forces was considered by Castanon \ea \cite{CASTANONHHO2020} while 
extension to non-Newtonian fluids was consedered by Botti \ea \cite{BottiMHHONN2021}.

In this work, we consider two HHO schemes that are novel variations of existing schemes with improved features.
In both cases, the Dirichlet condition on the velocity is enforced weakly in the spirit of \cite{Botti.Di-Pietro.ea:19*1}.
The first scheme borrows the skew-symmetric convective term treatment proposed by Botti \ea \cite{Botti.Di-Pietro.ea:19*1}, 
while the pressure-velocity coupling is inspired by the Hybridizable Discontinuous Galerkin (HDG) method of \cite{Rhebergen.Wells:18}.
Likewise the latter HDG formulation, the method is kinetic energy preserving and yields an exactly divergence free H-div conforming velocity approximation. 
Pressure-robustness is also guaranteed, meaning that the pressure error does not influence the velocity error.
We use polynomials of degree $k+1$ and $k$ to approximate the velocity and the pressure over mesh elements, respectively.
Polynomials of degree $k+1$, $k+1$ and $k$ are employed to approximate the numerical trace of the pressure on mesh faces, 
velocity on Neumann boundary faces and velocity on internal and Dirichlet boundary faces, respectively.
With respect to the HDG formulation of \cite{Rhebergen.Wells:18}, thanks to HHO treatment of the viscous terms in the spirit of \cite{Cockburn.Di-Pietro.ea:16} (see also \cite[Section 5.1]{Di-Pietro.Droniou:20}), 
improved orders of convergence are gained on simplicial meshes: in particular we observe $k+2$ and $k+1$ convergence rate for the velocity and pressure error in $L^2$ norm, respectively. 
Note that the leading block size of the statically condensed Jacobian matrix is 
driven by skeletal velocity DOFs, mostly associated to polynomials functions of degree $k$.
The second scheme targets robustness in the inviscid limit. 
Indeed, as we shall demonstrate by means of numerical test cases, the formulation is able to cope with the incompressible Euler equations.
Inspired by DG methods based on Godunov fluxes, see \eg \cite{Elsworth.Toro:1992,Bassi.Crivellini.ea:06,Bassi.Massa.ea:18,Massa.Ostrowski.ea:22} and following the ideas proposed by \cite{VilaPerez2021}
we employ an Harten, Lax and van Leer (HLL) approximated Riemann solver for designing pressure-velocity coupling and convective term treatment. 
The scheme relies on polynomials of degree $k$ for both elemental and skeletal unknowns resulting in convergence rates of order $k+1$ 
for both the velocity and pressure error in $L^2$ norm. Moreover, in the diffusion dominated regime, order $k+2$ for the velocity is recovered.

The paper is organized as follows. 
After introducing the discrete settings and the HHO discretization of the viscous term 
that the two formulations share in common, see Sections \ref{sec:meshSet}-\ref{sec:HHOLap}, 
the \hhohdiv and \hhohll formulations of the incompressible Navier-Stokes equations 
are provided in Sections \ref{sec:HHO-Hdiv} and \ref{sec:HHO-Hll}, respectively.
Notable features of the proposed formulations are outlined in Section \ref{subsec:hdivConf}, 
which focus on H-div conformity of {\tt HHO-Hdiv}, and Section \ref{subsec:hll_riemann},
which outlines the derivation of the HLL-type Riemann solver which \hhohll is based upon. 
In Section \ref{sec:numResults} we tackle several flow configurations, namely Kovasznay flow (Section \ref{sec:Kova}),
LLMS pressure gradient (Section \ref{sec:LLMS}), Gresho-Chan vortex (Section \ref{sec:GCVortex}),
double shear layer (\ref{sec:DSLayer}), and lid-driven cavity flow (Section {\ref{sec:LidDriven}).

\section{Two HHO methods for the Navier-Stokes problem}

\subsection{Discrete setting}
\label{sec:meshSet}

The HHO formulations proposed in this work allow to deal with two and three dimensional flow problems. 
Nevertheless, since numerical test cases focus on two space dimensions, discrete settings are provided for the 2D case.
We consider meshes of the domain $\Omega$ corresponding to couples $\Mh\coloneq(\Th,\Fh)$, where $\Th$ 
is a finite collection of polygonal elements such that $h\coloneq\max_{T\in\Th}h_T>0$ with $h_T$ denoting the diameter of $T$, 
while $\Fh$ is a finite collection of line segments. 
Extension to 3D requires to consider meshes composed of polyhedral elements and polygonal faces.
It is assumed henceforth that the mesh $\Mh$ is shape and contact regular, as detailed in \cite[Definition 1.4]{Di-Pietro.Droniou:20}.
For each mesh element $T \in \Th$, the faces contained in the element boundary $\partial T$ are collected in the set $\FT$,
and, for each mesh face $F \in \Fh$, $\TF$ is the set containing the one or two mesh elements sharing $F$.
We define three disjoint subsets of the set $\FT$:
the set of Dirichlet boundary faces $\FT^\dir \coloneq \{F \in \FT : F \subset \partial \Omega_\dir\}$;
the set of Neumann boundary faces $\FT^\neu \coloneq \{F \in \FT : F \subset \partial \Omega_\neu\}$;
the set of internal faces $\FT^\internal \coloneq \FT\setminus\big(\FT^\dir\cup\FT^\neu\big)$.
For future use, we also let $\FT^{\internal,\dir}\coloneq\FT^\internal\cup\FT^\dir$.
Using the same arguments we define three disjoint subsets of the set $\Fh$: $\Fh^\dir$, $\Fh^\neu$, $\Fh^\internal$
and we let $\Fh^{\internal,\dir}\coloneq\Fh^\internal\cup\Fh^\dir$.
For all $T\in\Th$ and all $F\in\FT$, $\normal_{TF}$ denotes the unit vector normal to $F$ pointing out of $T$.

Hybrid High-Order methods hinge on local polynomial spaces on mesh elements and faces.
For given integers $\ell\ge 0$ and $n\ge 1$, we denote by $\Polyd{n}{\ell}$ the space of $n$-variate polynomials of total degree $\le\ell$ (in short, of degree $\ell$).
For $X$ mesh element or face, we denote by $\Poly{\ell}(X)$ the space spanned by the restriction to $X$ of functions in $\mathbb{P}_d^\ell$.
When $X$ is a mesh face, the resulting space is isomorphic to $\mathbb{P}_{d-1}^\ell$ (see \cite[Proposition 1.23]{Di-Pietro.Droniou:20}).

Let again $X$ denote a mesh element or face.
The local $L^2$-orthogonal projector $\pi_X^\ell:L^2(X)\to \Poly{\ell}(X)$ is such that, for all $q\in L^2(X)$,
\[
\int_X ( q - \pi_X^\ell q ) r = 0\qquad\forall r\in\Poly{\ell}(X).
\]
Notice that, above and in what follows, we omit the measure from integrals as it can always be inferred from the context.
The $L^2$-orthogonal projector on $\Poly{\ell}(X)^d$, obtained applying $\pi_X^\ell$ component-wise, is denoted by $\boldsymbol{\pi}_X^\ell$.

\subsection{HHO discretization of the Laplace operator}
\label{sec:HHOLap}

The HHO discretizations of the Navier-Stokes problem considered in this work hinge on velocity reconstructions devised at the element level and obtained assembling diffusive potential reconstructions component-wise.
In what follows, we let a mesh element $T\in\Th$ be fixed, denote by $k\ge 0$ the degree of polynomials attached to internal and Dirichlet mesh faces, by $t \in\{k,k+1\}$ the degree of polynomials attached to mesh elements and by $f \in \{k,k+1\}$ the degree of polynomials attached to Neumann faces. 
We remark that, the motivation for choosing $ t $ and $ f $ one degree higher than $ k $ will be given in Section \ref{subsec:hdivConf}.

\subsubsection{Scalar potential reconstruction}

The velocity reconstruction is obtained leveraging, for each component, the \emph{scalar potential reconstruction} originally introduced in \cite{Di-Pietro.Ern.ea:14} in the context of scalar diffusion problems (see also \cite{Cockburn.Di-Pietro.ea:16} and \cite[Section 5.1]{Di-Pietro.Droniou:20} for its generalization to the case of different polynomial degrees on elements and faces).
Define the local scalar HHO space 
\begin{equation}\label{eq:VTk}
  \underline{V}_T^{t,k,f} \coloneq \left\{ \underline{v}_T = \left(v_T, (v_F)_{F \in \mathcal{F}_T} \right) : 
  \text{$v_T \in \Poly{t}(T)$, $v_F \in \Poly{k}(F)$ for all $F \in\mathcal{F}^{\internal,\dir}_T$ and $v_F \in \Poly{f}(F)$ for all $F \in\mathcal{F}^{\neu}_T$}
  \right\}.  
\end{equation}
The scalar potential reconstruction operator $\mathfrak{p}_T^{k+1}$: $\underline{V}_T^{t,k,f} \rightarrow \Poly{k+1}(T)$ maps a vector of polynomials of $\underline{V}_T^{t,k,f}$ onto a polynomial of degree $(k+1)$ over $T$ as follows:
Given $\underline v_{T} \in \underline{V}_T^{t,k,f}$, $\mathfrak{p}_T^{k+1}\underline{v}_T$ is the unique polynomial in $\Poly{k+1}(T)$ satisfying
\begin{equation*} 
  \begin{alignedat}{2}
    \int_T \nabla \mathfrak{p}_T^{k+1} \underline{v}_T \cdot \nabla {w_T}
    &=  \int_T \nabla {v_T} \cdot \nabla {w_T} + \sum_{F \in \FT} \int_F \left({v_F} - {v_T}\right) \, \nabla {w_T} \cdot \normal_{TF}
    &\qquad&  \forall w_T \in \Poly{k +1}(T),
    \\
    \int_T \mathfrak{p}_T^{k+1} \underline{v}_T &= \int_T {v_T}.
  \end{alignedat}
\end{equation*}
Computing $\mathfrak{p}_T^{k+1}$ for each $T\in\Th$ requires to solve a small linear system.
This task shows optimal scalability on parallel architectures.

\subsubsection{Velocity reconstruction}
\label{sec::HHOVelRec}

Define, in analogy with \eqref{eq:VTk}, the following vector-valued HHO space for the velocity:
\begin{equation}\label{eq:VTkV}
  \underline{\boldsymbol{V}}_T^{t,k,f} \coloneq \left\{ \underline{\vVec}_T = \big(\vVec_T, (\vVec_F)_{F \in \mathcal{F}_T} \big) : 
  \text{%
    ${\vVec}_T \in \Poly{t}(T)^d$, $\vVec_F  \in \Poly{k}(F)^d$ for all $F \in\mathcal{F}_T^{\internal,\dir}$ and $\vVec_F  \in \Poly{f}(F)^d$ for all $F \in\mathcal{F}_T^{\neu}$
  }
  \right\}.
\end{equation}
The \emph{velocity reconstruction} $\mathfrak{P}_T^{k+1}$: $\underline{\boldsymbol{V}}_T^{t,k,f} \rightarrow \Poly{k+1}(T)^d$ is obtained setting
\begin{equation*} 
  \mathfrak{P}_T^{k+1} \underline{\vVec}_T \coloneq \big(\mathfrak{p}_T^{k+1} \underline{v}_{T,i}\big)_{i=1,\ldots,d},
\end{equation*}
where, for all $i=1,\ldots,d$, $\underline{v}_{T,i}\in\underline{V}_T^{t,k,f}$ is obtained gathering the $i$th components of the polynomials 
in $\underline{\vVec}_T$, \ie $\underline{v}_{T,i}\coloneq \big(v_{T,i}, (v_{F,i})_{F\in\FT}\big)$ if $\vVec_T=(v_{T,i})_{i=1,\ldots,d}$ and $\vVec_F=(v_{F,i})_{i=1,\ldots,d}$ for all $F\in\FT$.

 The local interpolation operator $\underline{\boldsymbol{I}}_T^{t,k,f}: H^1(T)^d \rightarrow \underline{\boldsymbol{V}}_T^{t,k,f}$ is defined as follows:
 For all $\vVec \in H^1(T)^d$,
 \begin{equation} 
   \label{def:intOpHHO}
   \underline{\boldsymbol{I}}_T^{t,k,f} \vVec \coloneq \big( \boldsymbol{\pi}^{t}_T \vVec, (\boldsymbol{\pi}^k_F \vVec_{|F})_{F \in \FT^{\internal,\dir}}, (\boldsymbol{\pi}^{f}_F \vVec_{|F})_{F \in \FT^{\neu}}\big).
 \end{equation}
 Following \cite[Section 5.1.3]{Di-Pietro.Droniou:20}, it is possible to demonstrate that the velocity reconstruction is such that,
 for all $\vVec \in H^1(T)^d$
 \begin{equation} 
   \text{%
     $\int_T \big(\nabla \mathfrak{P}_T^{k+1} \underline{\boldsymbol{I}}_T^{t,k,f} \vVec - \nabla \vVec\big) : \nabla \wVec_T = 0$ for all $\wVec_T \in \Poly{k+1}(T)^d$
     and $\int_T\mathfrak{P}_T^{k+1}\underline{\boldsymbol{I}}_T^{t,k,f}\vVec = \int_T\vVec$.
   }
 \end{equation}
 The above result is of crucial importance for inferring the approximation properties of the velocity reconstruction operator.

\subsubsection{Face residuals}

Let $T\in\Th$ and $F\in\FT$.
The stabilization bilinear form for the HHO discretization of the viscous term in the momentum equation \eqref{stokesProb:momentum} hinges 
on the residual operator $\mathfrak{r}_{TF}^{k}: \underline{V}_T^{t,k,f} \rightarrow \underline{V}_T^{t,k,f}$ 
defined as follows 
\[
\mathfrak{r}_{TF}^{k} \underline{v}_T \coloneq  \mathfrak{r}_{F}^{k,f} \underline{v}_T - \mathfrak{r}_{F}^{t} \underline{v}_T 
\]
where the \emph{face residual} 
$\mathfrak{r}_{F}^{k,f} : \underline{V}_T^{t,k,f} \rightarrow \Poly{k,f}(F)$ 
and the \emph{element residual}
$\mathfrak{r}_{T}^{t} : \underline{V}_T^{t,k,f} \rightarrow \Poly{t}(T)$ 
such that, for all $\underline{v}_T\in\underline{V}_T^{t,k,f}$,
\[
\mathfrak{r}_{F}^{k,f} \underline{v}_T \coloneq \begin{cases} 
\pi^{k}_F \big(v_F - \mathfrak{p}_T^{k+1} \underline{v}_T \big) & \text{if $F \in \Fh^{\internal,\dir}$} \\ 
\pi^{f}_F \big(v_F - \mathfrak{p}_T^{k+1} \underline{v}_T \big) & \text{if $F \in \Fh^{\neu}$}
\end{cases} \quad
\text{and} \quad
\mathfrak{r}_{T}^{t} \underline{v}_T = 
\coloneq \pi^{t}_T \big(v_T - \mathfrak{p}_T^{k+1} \underline{v}_T \big).
\]

The vector residual $\mathfrak{R}_{TF}^{k}: \underline{\boldsymbol{V}}_T^{t,k,f} \rightarrow \underline{\boldsymbol{V}}_T^{t,k,f}$ 
is such that, for all $\underline{\vVec}_T\in\underline{\boldsymbol{V}}_T^{t,k,f}$:
\[
\mathfrak{R}_{TF}^{k} \underline{\vVec}_T \coloneq \big( \mathfrak{r}_{TF}^{k} \underline{v}_{T,i} \big)_{i=1,\ldots,d}.
\]

\subsection{Local and global HHO spaces for velocity and pressure unknowns}
In addition to the scalar and vector valued local velocity spaces $\underline{V}_T^{t,k,f}$ and $\underline{\boldsymbol{V}}_T^{t,k,f}$ defined in \eqref{eq:VTk} and \eqref{eq:VTkV}, respectively,
we need to introduce local pressure spaces and global HHO spaces for both velocity and pressure.

Let $T\in\Th$ and define the local scalar HHO space for the pressure unknown
\begin{equation}\label{eq:QTk}
  \underline{Q}_T^{k,p} \coloneq \left\{ \underline{q}_T = \left(q_T, (q_F)_{F \in \mathcal{F}_T} \right) : 
  \text{$q_T \in \Poly{k}(T)$ and $q_F \in \Poly{p}(F)$ for all $F \in\mathcal{F}_T$}
  \right\},  
\end{equation}
where $p\in\{k,k+1\}$.
The global scalar HHO space is defined and follows
\[
\underline{Q}_h^{k,p}\coloneq\left\{
\underline{q}_h=\big((q_T)_{T\in\Th}, (q_F)_{F\in\Fh}\big)\,:\,
\text{$q_T\in\Poly{k}(T)$ for all $T\in\Th$ and $q_F\in\Poly{p}(F)$ for all $F\in\Fh$}
\right\}.
\]
For all $\underline{q}_h\in\underline{\boldsymbol{Q}}_h^{k,p}$ and all $T\in\Th$, we denote by $\underline{q}_T\in\underline{\boldsymbol{Q}}_T^{k,p}$ the restriction of $\underline{q}_h$ to $T$.

To conclude, the global vector HHO space reads
\[
\underline{\boldsymbol{V}}_h^{t,k,f}\coloneq \left\{ \underline{\vVec}_h = \big( (\vVec_T)_{T\in\Th} , (\vVec_F)_{F \in \Fh} \big) : 
\begin{array}{l}
\text{${\vVec}_T \in \Poly{t}(T)^d$ for all $T\in\Th$},\\ \text{$\vVec_F  \in \Poly{k}(F)^d$ for all $F \in \Fh^{\internal,\dir}$,} \\ \text{$\vVec_F  \in \Poly{f}(F)^d$ for all $F \in \Fh^{\neu}$
}
\end{array}
\right\}.
\]
For all $\underline{\vVec}_h\in\underline{\boldsymbol{V}}_h^{t,k,f}$ and all $T\in\Th$, we denote by $\underline{\vVec}_T\in\underline{\boldsymbol{V}}_T^{t,k,f}$ the restriction of $\underline{\vVec}_h$ to $T$.

\subsection{A pointwise divergence free H-div conforming HHO scheme}
\label{sec:HHO-Hdiv}
\subsubsection{Local and global residuals}
We combine the HHO discretization of the viscous term with a hybrid approximation of the pressure inspired by \cite{Rhebergen.Wells:18}.
Given $(\underline \uVec_T, \underline{p}_T)\in \underline{\boldsymbol{V}}_T^{t,k,f}\times \underline{Q}_T^{k,p}$, the local residuals
$r^\momentum_{I,T} ((\underline{\boldsymbol{u}}_T,\underline{p}_T);\cdot):\underline{\boldsymbol{V}}_T^{t,k,f}\to\Real$ of the space discrete momentum
and $ r^\mass_{I,T}(\underline{\boldsymbol{u}}_T;\cdot):\underline{Q}^{k,p}_T\to\Real$ of the discrete mass conservation equations are such that, 
for all $\underline \vVec_T\in\underline{\boldsymbol{V}}_T^{t,k,f}$ and all $\underline{q}_T\in\underline{Q}_T^{k,p}$,
\begin{subequations}
\begin{align*} 
  r^\momentum_{I,T} ((\underline{\boldsymbol{u}}_T,\underline{p}_T);\underline{\boldsymbol{v}}_T)
  &\coloneq
  \begin{aligned}[t]
    &\int_T \nu \nabla \mathfrak{P}_T^{k+1} \underline{\uVec}_{T} : \nabla \mathfrak{P}_T^{k+1} \underline{\vVec}_{T}
    + \sum_{F \in \FT} \frac{\nu}{h_F} \int_F \mathfrak{R}_{TF}^k \underline{\vVec}_T \, \cdot \, \mathfrak{R}_{TF}^k \underline{\vVec}_T
    \\
        & - \sum_{F \in \FT^\dir} \int_F \left[
      \bigl( \normal_{TF} \cdot \nu \nabla \mathfrak{P}_T^{k+1} \underline{\uVec}_{T} \bigr) \cdot \vVec_F   
      + \uVec_F \cdot \bigl( \normal_{TF} \cdot \nu \nabla \mathfrak{P}_T^{k+1} \underline{\vVec}_{T}\bigr)
      \right]
    +\sum_{F \in \FT^\dir} \frac{\eta \nu}{h_F} \int_F \uVec_F \cdot \vVec_{F}   
    \\
    &
    - \int_T p_{T} \, (\nabla \cdot \vVec_T)
    + 
      \sum_{F \in \FT} \int_F p_F \, (\vVec_T - \vVec_F) \cdot \normal_{TF}
    + \sum_{F \in \FT^\dir} \int_F p_F \, (\vVec_F \cdot \normal_{TF})
\\
&
 + \frac{1}{2} \left( \int_T (\uVec_T \cdot \nabla \uVec_T) \cdot \vVec_T - \int_T (\uVec_T \otimes \uVec_T) : \nabla \vVec_T  + \sum_{F \in \FT} \int_F (\uVec_F \cdot \normal_{TF}) \, ( \uVec_F \cdot \vVec_T - \uVec_T \cdot \vVec_F) \right)
\\
&
 + \frac{1}{2} \left( \sum_{F \in \FT^\dir} \int_F (\uVec_F \cdot \normal_{TF}) \, \boldsymbol{g}_\dir  \cdot \vVec_F  + \sum_{F \in \FT^\neu} \int_F (\uVec_F \cdot \normal_{TF}) \, \uVec_T \cdot \vVec_F) \right)
\\
&
    - \sum_{F \in \FT^\dir} \int_F
    \boldsymbol{g}_\dir\cdot \left(
      \normal_{TF} \cdot \nu \nabla \mathfrak{P}_T^{k+1} \underline{\vVec}_{T}
      + \frac{\eta \nu}{h_F} \vVec_{F}
      \right)
- \sum_{F \in \FT^\neu} \int_F \boldsymbol{g}_\neu \cdot \vVec_F
    + \int_{\aT} \left(\frac{\partial \uVec_{T}}{\partial t}-\boldsymbol{f}\right) \cdot \vVec_T,
  \end{aligned}
  \\ 
  r^\mass_{I,T}(\underline{\boldsymbol{u}}_T;\underline{q}_T)
  &\coloneq
    \begin{aligned}[t]
  & - \int_T  (\nabla\cdot{\uVec}_T) \, q_T
  +   
    \sum_{F \in \FT} \int_F ( \uVec_T - \uVec_F ) \cdot \normal_{TF} \, q_F 
      + \sum_{F \in \FT^\dir} \int_F \left(\uVec_F - \boldsymbol{g}_\dir \right) \cdot \normal_{TF}\, q_F
    \end{aligned}
\end{align*}
\end{subequations}
In the expression of $r^\momentum_{I,T}((\underline{\boldsymbol{u}}_T,\underline{p}_T);\cdot)$, $\eta>0$ is a user-dependent parameter that has to be taken large enough to ensure coercivity.
The penalty term where the parameter $\eta$ appears, along with the consistency terms in the second line and the terms involving the boundary datum $\boldsymbol{g}_\dir$ in the fifth and sixth lines, 
are responsible for the weak enforcement of Dirichlet boundary conditions for the velocity.
In the numerical tests provided below, $\eta$ is taken equal to $3$.

The global residuals
$r^\momentum_{I,h}((\underline{\boldsymbol{u}}_h,\underline{p}_h);\cdot):\underline{\boldsymbol{V}}^{t,k,f}_h\to\Real$ and
$r^\mass_{I,h}(\underline{\boldsymbol{u}}_h;\cdot):\underline{Q}^{k,p}_h\to\Real$ are obtained by element-by-element assembly of the local residuals, \textit{i.e.}:
For all $\underline{\boldsymbol{v}}_h\in\underline{\boldsymbol{V}}_h^{t,k,f}$ and all $\underline{q}_h\in\underline{Q}_h^{k,p}$,
\begin{equation}
\label{HHOdisrcGlobalResHDIV}
r_{I,h}^\momentum\left((\underline{\boldsymbol{u}}_h,\underline{p}_h);\underline{\boldsymbol{v}}_h\right)
\coloneq \sum_{T\in\Th} r^\momentum_{I,T}\left((\underline{\boldsymbol{u}}_T,\underline{p}_T);\underline{\boldsymbol{v}}_T\right),\qquad
r_{I,h}^\mass(\underline{\boldsymbol{u}}_h;\underline{q}_h) \coloneq \sum_{T\in\Th} r^\mass_{I,T}(\underline{\boldsymbol{u}}_T;\underline{q}_{T}).
\end{equation}

The HHO scheme with pointwise divergence free H-div conforming velocity field is obtained setting the polynomial degrees as $t=f=p=k+1$.
\begin{scheme}[{\tt HHO-Hdiv}] \label{scheme:hho.hdiv}
  Find $(\underline{\boldsymbol{u}}_h,\underline{p}_h)\in\underline{\boldsymbol{V}}_h^{k+1,k,k+1}\times\underline{Q}^{k,k+1}_h$ such that
  \begin{equation}\label{eq:hho.hdiv}
    \begin{alignedat}{2}
      r_{I,h}^\momentum((\underline{\boldsymbol{u}}_h,\underline{p}_h);\underline{\boldsymbol{v}}_h) &= 0 &\qquad& \forall \underline{\boldsymbol{v}}_h\in\underline{\boldsymbol{V}}_h^{k+1,k,k+1},
      \\
      r_{I,h}^\mass(\underline{\boldsymbol{u}}_h;\underline{q}_h) &= 0 &\qquad& \forall \underline{q}_h\in\underline{Q}_h^{k,k+1}.
    \end{alignedat}
  \end{equation}
\end{scheme}

\subsubsection{Discrete mass conservation and H-div conformity}
\label{subsec:hdivConf}
For each $\underline{\uVec}_h \in \underline{\boldsymbol{V}}_h^{t,k,f}$, $r_{I,h}^{\mass} = 0$ for all $\underline{q}_h \in \underline{Q}_h^{k,p}$
implies that
\begin{enumerate}
\item
if $q_F = 0$ for all $F \in \Fh$, then $\int_T  (\nabla\cdot{\uVec}_T) \, q_T = 0, \; \text{for each $T \in \Th$}$;
\item
if $q_T = 0$ for all $T \in \Th$, then $ \left\{ \begin{aligned} 
\int_F  (\uVec_{T} - \uVec_{T'}) \cdot \normal_{TF} \, q_F = 0, \; &\text{ for each $F\in\FT^\internal\cap\mathcal{F}_{T'}^\internal$ with $T,T'\in\Th$, $T\neq T'$};\\  
\int_F  (\uVec_T - \uVec_F) \cdot \normal_{TF} \, q_F = 0, \; &\text{ for each $T \in \Th, F \in \FT^\neu$};\\
\int_F  (\uVec_T - \boldsymbol{g}_\dir) \cdot \normal_{TF} \, q_F = 0, \; &\text{ for each $T \in \Th, F \in \FT^\dir$}.
\end{aligned} \right.$
\end{enumerate}
Setting $t=f=p=k+1$, it is straightforward to conclude that 
\begin{align}
\nabla \cdot \uVec_T &= 0, &\forall \xVec \in T, \quad&\forall T \in \Th; & \label{hdivCondVol}\\
(\uVec_{T} - \uVec_{T'}) \cdot \normal_{TF} &= 0, &\forall \xVec \in F, \quad&\text{ for each $F\in\FT^\internal\cap\mathcal{F}_{T'}^\internal$ with $T,T'\in\Th$, $T\neq T'$}; & \label{eq:hdivCondIfaces}\\
(\uVec_{T} - \uVec_{F}) \cdot \normal_{TF} &= 0, &\forall \xVec \in F,\quad& \text{ for each $T \in \Th, F \in \FT^\neu$}; &\label{eq:hdivCondNfaces}\\
(\uVec_{T} - \boldsymbol{\pi}_F^{k+1}{\boldsymbol{g}_\dir}) \cdot \normal_{TF} &= 0, &\forall \xVec \in F, \quad&\text{ for each $T \in \Th, F \in \FT^\dir$}. \label{eq:hdivCondDfaces}&
\end{align}

\begin{remark}[Polynomial degree of velocity unknowns over Neumann faces]
When $t=p=k+1$, choosing $f=k+1$ is mandatory to ensure H-div conformity on Neumann faces without constraining velocity unknowns over mesh elements. 
Note in particular that setting $f=k$ would imply that, according to \eqref{eq:hdivCondNfaces}, $\uVec_{T}|_F \in \Poly{k}(F)$ for each $F \in \Fh^\neu$.
We also remark that setting $f=k$ is perfectly fine if $\boldsymbol{g}_\neu$ is a constant, in particular traction-free boundary conditions commonly employed 
on outflow sections of channel flows can be dealt with without raising the polynomial degree on outflow faces. 
\end{remark}

\begin{remark}[$h$-Convergence rates on simplicial and general meshes]
The \hhohdiv formulations show convergence rates of order $k+1$ for the pressure error in $L^2$ norm on both simplicial and general meshes. 
Convergence rates of order $k+2$ and $k+1$ for the velocity and velocity gradients error in $L^2$ norm, respectively, are attained on simplicial meshes, while, on general meshes, convergence rates are reduced by one order.
Choosing $t=k+1$ and $f=p=k$ allows to recover optimal convergence rates for the velocity on general meshes but H-div conformity and pressure-robustness are lost. 
Indeed, \eqref{eq:hdivCondIfaces}-\eqref{eq:hdivCondNfaces}-\eqref{eq:hdivCondDfaces} are replaced by the following looser conditions on the normal trace of the velocity 
\begin{align}
(\boldsymbol{\pi}_F^{k} \uVec_{T} - \boldsymbol{\pi}_F^{k}\uVec_{T'}) \cdot \normal_{TF} &= 0, &\forall \xVec \in F, \quad&\text{ for each $F\in\FT^\internal\cap\mathcal{F}_{T'}^\internal$ with $T,T'\in\Th$, $T\neq T'$}; & \label{eq:hdivCondIfacesK}\\
(\boldsymbol{\pi}_F^{k} \uVec_{T} - \uVec_{F}) \cdot \normal_{TF} &= 0, &\forall \xVec \in F,\quad& \text{ for each $T \in \Th, F \in \FT^\neu$} &\label{eq:hdivCondNfacesK}\\
(\boldsymbol{\pi}_F^{k} \uVec_{T} - \boldsymbol{\pi}_F^{k}{\boldsymbol{g}_\dir}) \cdot \normal_{TF} &= 0, &\forall \xVec \in F, \quad&\text{ for each $T \in \Th, F \in \FT^\dir$} \label{eq:hdivCondDfacesK}&
\end{align}
Since \eqref{hdivCondVol} holds, the method yields a pointwise divergence free velocity field. 
The HHO scheme with pointwise divergence free velocity and optimal convergence properties on general meshes is outlined in what follows for the sake of completeness.
\begin{scheme}[{\tt HHO-DivFree}] \label{scheme:hho.divfree}
  Find $(\underline{\boldsymbol{u}}_h,\underline{p}_h)\in\underline{\boldsymbol{V}}_h^{k+1,k,k}\times\underline{Q}^{k,k}_h$ such that
  \begin{equation}\label{eq:hho.divfree}
    \begin{alignedat}{2}
      r_{I,h}^\momentum((\underline{\boldsymbol{u}}_h,\underline{p}_h);\underline{\boldsymbol{v}}_h) &= 0 &\qquad& \forall \underline{\boldsymbol{v}}_h\in\underline{\boldsymbol{V}}_h^{k+1,k,k+1},
      \\
      r_{I,h}^\mass(\underline{\boldsymbol{u}}_h;\underline{q}_h) &= 0 &\qquad& \forall \underline{q}_h\in\underline{Q}_h^{k,k+1}.
    \end{alignedat}
  \end{equation}
\end{scheme}
\end{remark}

\subsection{A HHO formulation robust in the inviscid limit}
\label{sec:HHO-Hll}

\subsubsection{Local and global residuals}

Given $(\underline \uVec_T, \underline{p}_T)\in \underline{\boldsymbol{V}}_T^{t,k,f}\times \underline{Q}_T^{k,p}$, the local residuals
$r^\momentum_{II,T} ((\underline{\boldsymbol{u}}_T,\underline{p}_T);\cdot):\underline{\boldsymbol{V}}_T^{t,k,f}\to\Real$ of the space discrete momentum
and $ r^\mass_{II,T}(\underline{\boldsymbol{u}}_T;\cdot):\underline{Q}^{k,p}_T\to\Real$ of the discrete mass conservation equations are such that, 
for all $\underline \vVec_T\in\underline{\boldsymbol{V}}_T^{t,k,f}$ and all $\underline{q}_T\in\underline{Q}_T^{k,p}$,
\begin{subequations}
\begin{align*} 
  r^\momentum_{II,T} ((\underline{\boldsymbol{u}}_T,\underline{p}_T);\underline{\boldsymbol{v}}_T)
  &\coloneq
  \begin{aligned}[t]
    &\int_T \nu \nabla \mathfrak{P}_T^{k+1} \underline{\uVec}_{T} : \nabla \mathfrak{P}_T^{k+1} \underline{\vVec}_{T}
    + \sum_{F \in \FT} \frac{\nu}{h_F} \int_F \mathfrak{R}_{TF}^k \underline{\vVec}_T \, \cdot \, \mathfrak{R}_{TF}^k \underline{\vVec}_T
    \\
        & - \sum_{F \in \FT^\dir} \int_F \left[
      \bigl( \normal_{TF} \cdot \nu \nabla \mathfrak{P}_T^{k+1} \underline{\uVec}_{T} \bigr) \cdot \vVec_F   
      + \uVec_F \cdot \bigl( \normal_{TF} \cdot \nu \nabla \mathfrak{P}_T^{k+1} \underline{\vVec}_{T}\bigr)
      \right]
    +\sum_{F \in \FT^\dir} \frac{\eta \nu}{h_F} \int_F \uVec_F \cdot \vVec_{F}   
    \\
&
    - \sum_{F \in \FT^\dir} \int_F
    \boldsymbol{g}_\dir\cdot \left(
      \normal_{TF} \cdot \nu \nabla \mathfrak{P}_T^{k+1} \underline{\vVec}_{T}
      + \frac{\eta \nu}{h_F} \vVec_{F}
      \right)
    - \int_T p_{T} \, (\nabla \cdot \vVec_T)
    + \sum_{F \in \FT} \int_F p_F \, (\vVec_T \cdot \normal_{TF})
    \\
    &
 - \int_T (\uVec_T \otimes \uVec_T) : \nabla \vVec_T  + \sum_{F \in \FT} \int_F (\uVec_F \cdot \normal_{TF}) \, \uVec_F \cdot \vVec_T
    + \sum_{F \in \FT} \int_F s^+ (\uVec_T - \uVec_F) \cdot (\vVec_T - \vVec_F) 
\\
&
    + \sum_{F \in \FT^\dir} \int_F s^- (\boldsymbol{g}_\dir - \uVec_F) \cdot \vVec_F 
- \sum_{F \in \FT^\neu} \int_F (\boldsymbol{g}_\neu - p_F \normal_{TF}) \cdot \vVec_F
    + \int_{\aT} \left(\frac{\partial \uVec_{T}}{\partial t}-\boldsymbol{f}\right) \cdot \vVec_T,
  \end{aligned}
  \\ 
  r^\mass_{II,T}(\underline{\boldsymbol{u}}_T;\underline{q}_T)
  &\coloneq
    \begin{aligned}[t]
  & - \int_T  {\uVec}_T \cdot \nabla q_T
  +   
    \sum_{F \in \FT} \int_F \uVec_F \cdot \normal_{TF} \, q_T 
    + \sum_{F \in \FT} \int_F \dfrac{s^+}{a^2}  (p_T - p_F)  \, (q_T - q_F) 
      \\
   & + \sum_{F \in \FT^\neu} \int_F \dfrac{s^-}{a^2} \left( ( \boldsymbol{g}_\neu + \normal_{TF} \cdot \nu \nabla \mathfrak{P}_T^{k+1} \underline{\uVec}_{T})  \cdot \normal_{TF} - p_F \right) \, q_F,
    \end{aligned}
\end{align*}
\end{subequations}
where 
\begin{equation} 
\label{eq:waveSpeed}
s^\pm = \frac{1}{2} \left( \uVec_F \cdot \normal_{TF} \pm \sqrt{(\uVec_F \cdot \normal_{TF})^2 + 4 a^2} \right)
\end{equation} 
and $a$ is the artificial compressibility parameter.
In the numerical tests provided below, $a$ is taken as a unit velocity.
In the expression of $r^\momentum_{II,T}((\underline{\boldsymbol{u}}_T,\underline{p}_T);\cdot)$, $\eta>0$ is a user-dependent parameter that has to be taken large enough to ensure coercivity.
In the numerical tests provided below, $\eta$ is taken equal to $3$.

The global residuals
$r^\momentum_{II,h}((\underline{\boldsymbol{u}}_h,\underline{p}_h);\cdot):\underline{\boldsymbol{V}}^{t,k,f}_h\to\Real$ and
$r^\mass_{II,h}(\underline{\boldsymbol{u}}_h;\cdot):\underline{Q}^{k,p}_h\to\Real$ are obtained by element-by-element assembly of the local residuals, \textit{i.e.}:
For all $\underline{\boldsymbol{v}}_h\in\underline{\boldsymbol{V}}_h^{t,k,f}$ and all $\underline{q}_h\in\underline{Q}_h^{k,p}$,
\begin{equation}
\label{HHOdisrcGlobalResHLL}
r_{II,h}^\momentum\left((\underline{\boldsymbol{u}}_h,\underline{p}_h);\underline{\boldsymbol{v}}_h\right)
\coloneq \sum_{T\in\Th} r^\momentum_{II,T}\left((\underline{\boldsymbol{u}}_T,\underline{p}_T);\underline{\boldsymbol{v}}_T\right),\qquad
r_{II,h}^\mass(\underline{\boldsymbol{u}}_h;\underline{q}_h) \coloneq \sum_{T\in\Th} r^\mass_{II,T}(\underline{\boldsymbol{u}}_T;\underline{q}_{T}).
\end{equation}

The HHO scheme based on the HLL Riemann solver is obtained setting the polynomial degrees as $t=f=p=k$.
\begin{scheme}[{\tt HHO-HLL}] \label{scheme:hho.hll}
  Find $(\underline{\boldsymbol{u}}_h,\underline{p}_h)\in\underline{\boldsymbol{V}}_h^{k,k,k}\times\underline{Q}^{k,k}_h$ such that
  \begin{equation}\label{eq:hho.hll}
    \begin{alignedat}{2}
      r_{II,h}^\momentum((\underline{\boldsymbol{u}}_h,\underline{p}_h);\underline{\boldsymbol{v}}_h) &= 0 &\qquad& \forall \underline{\boldsymbol{v}}_h\in\underline{\boldsymbol{V}}_h^{k,k,k},
      \\
      r_{II,h}^\mass(\underline{\boldsymbol{u}}_h;\underline{q}_h) &= 0 &\qquad& \forall \underline{q}_h\in\underline{Q}_h^{k,k}.
    \end{alignedat}
  \end{equation}
\end{scheme}

\subsubsection{Artificial compressibility based HLL-type numerical flux}
\label{subsec:hll_riemann}

In order to outline the convective term and pressure-velocity coupling formulations we focus on the inviscid framework 
and, for the sake of conciseness, we assume that periodic boundary conditions are imposed over $\partial \Omega$.
Accordingly, the local residuals for the \hhohll discretization of the Euler equations on a periodic domain reads:
given $(\underline \uVec_T, \underline{p}_T)\in \underline{\boldsymbol{V}}_T^{k,k,k}\times \underline{Q}_T^{k,k}$,
for all $\underline \vVec_T\in\underline{\boldsymbol{V}}_T^{k,k,k}$ and all $\underline{q}_T\in\underline{Q}_T^{k,k}$
\begin{align} 
  r^{\momentum}_{II,T}((\underline{\boldsymbol{u}}_T,\underline{p}_T);\underline{\boldsymbol{v}}_T)
	& =
	\int_T \left(\uVec_{T}\otimes\uVec_{T}  +p_T \I \right): 
	\nabla \boldsymbol{v}_{T} + \int_{\aT} \left( \dfrac{\partial \uVec_T}{\partial t} - \boldsymbol{f} \right) \cdot \vVec_T 
	\nonumber \\
        & +	
	\sum_{F \in \FT}\int_F \left[\left(\uVec_{F}\otimes\uVec_{F}  +p_F \I\right)\cdot \normal_{TF}
 +s^+ ( \uVec_T - \uVec_F )\right] \cdot 
	(\boldsymbol{v}_{T} - \boldsymbol{v}_{F}),
	\label{eq:HLL_numerical_flux_qdm}\\
  r^{\mass}_{II,T}(\underline{\boldsymbol{u}}_T;\underline{q}_T)
	& =
  - \int_T  {\uVec}_T \cdot \nabla q_T +
	\sum_{F \in \FT}\int_F 
	\left[\uVec_F \cdot \normal_{TF} +\dfrac{s^+}{a^2} (p_T - p_F) \right] 
	\,(q_T - q_F).\label{eq:HLL_numerical_flux_cnt}
\end{align}
The terms within square brackets are the, so called, \emph{numerical fluxes} of the HHO formulation. 

In order to outline the trace conditions enforced over mesh faces we now consider the global residuals of the Euler equations 
with Dirichlet and Neumann boundary conditions imposed over $\partial \Omega$. 
For each $(\underline \uVec_h, \underline{p}_h)\in \underline{\boldsymbol{V}}_h^{k,k,k}\times \underline{Q}_h^{k,k}$,
$r_{II,h}^{\momentum} = 0$ for all $\underline{\vVec}_h \in \underline{V}_h^{k,k,k}$ such that $\vVec_T = 0$ for all $T \in \Th$,
implies that
\begin{equation*}
	\begin{aligned}
	&\int_F
	\left[s^+ ( \uVec_T - \uVec_F )  -s^- ( \uVec_{T'} - \uVec_F ) \right]\cdot 
	\boldsymbol{v}_F = 0 ,&
	 \; &\text{ for each $F\in\FT^\internal\cap\mathcal{F}_{T'}^\internal$ with $T,T'\in\Th$, $T\neq T'$};&\\
	&\int_F
	s^+ ( \uVec_T - \uVec_F )\cdot 
	\boldsymbol{v}_F = 0 ,&
	 \; &\text{ for each $T\in\Th$, $F\in\mathcal{F}_{T}^\neu$};&\\
	&\int_F
	\left[s^+ ( \uVec_T - \uVec_F )  -s^- ( \boldsymbol{g}_\dir - \uVec_F ) \right]\cdot 
	\boldsymbol{v}_F = 0 ,&
	 \; &\text{ for each $T\in\Th$, $F\in\mathcal{F}_{T}^\dir$}.&
	 \end{aligned}
\end{equation*}
For each $\underline \uVec_h \in \underline{\boldsymbol{V}}_h^{k,k,k}$,
$r_{II,h}^{\mass} = 0$ for all $\underline{q}_h \in \underline{Q}_h^{k,k}$ such that $q_T = 0$ for all $T \in \Th$,
implies that
\begin{equation*}
	\begin{aligned}
	&\int_F
	\left[s^+ ( p_T - p_F )  -s^- ( p_{T'} - p_F ) \right] 
	\,q_F = 0 ,&
	 \; &\text{ for each $F\in\FT^\internal\cap\mathcal{F}_{T'}^\internal$ with $T,T'\in\Th$, $T\neq T'$};&\\
	&\int_F
	\left[s^+ ( p_T - p_F )  -s^- \left( \boldsymbol{g}_\neu \cdot \normal_{TF} - p_F \right) \right] 
	\,q_F = 0 ,&
	 \; &\text{ for each $T\in\Th$, $F\in\mathcal{F}_{T}^\neu$};&\\
	&\int_F
	s^+( p_T - p_F ) 
	\,q_F = 0 ,&
	 \; &\text{ for each $T\in\Th$, $F\in\mathcal{F}_{T}^\dir$}.&  
	 \end{aligned}
\end{equation*}
Accordingly, the intermediate state at each $F\in\FT^\internal\cap\mathcal{F}_{T'}^\internal$ 
is the weighted average between the two neighbour element ($T,T'$) states, \ie
\begin{equation}\label{eq:face_variables}
	\boldsymbol{\pi}_F^{k}\left((s^{+}-s^{-})\uVec_{F}\right)=\boldsymbol{\pi}_F^{k}\left(s^{+}\uVec_{T}  -s^{-}\uVec_{T'}\right)
	\qquad
	\text{and}
	\qquad
	\pi_F^{k}\left((s^{+}-s^{-})p_{F}\right)=\pi_F^{k}\left(s^{+}p_{T} -s^{-}p_{T'}\right).
\end{equation}

The \hhohll formulation falls into the generalized framework of Riemann solvers for hybrid DG methods proposed by Vila-P\'erez \ea \cite{VilaPerez2021}
and is based on the Exact Riemann Solver (ERS) for variable density incompressible flows devised in \cite{Bassi.Massa.ea:18}.
In particular, following \cite{Toro:2009}, the ERS based HLL Riemann solver can be derived:
the flux$_\star$ in the \emph{star} region delimited by the two external acoustic waves reads
\begin{align*}
	\left(\uVec\otimes\uVec  +p\I \right)_{\star} \cdot\normal_{FT}
		&=	\dfrac{s_{T'}\left(\uVec_{T}\otimes\uVec_{T}  +p_{T}\I \right)
			-s_{T}\left(\uVec_{T'}\otimes\uVec_{T'}  +p_{T'}\I \right)}{s_{T'}-s_{T}} \cdot\normal_{FT}
			+\dfrac{s_{T'}s_{T}}{s_{T'}-s_{T}}\left(\uVec_{T'}  -\uVec_{T}\right),\\
	\left(\uVec\right)_{\star} \cdot\normal_{FT}
		&=	\dfrac{s_{T'}\uVec_{T}
			-s_{T}\uVec_{T'}}{s_{T'}-s_{T}} \cdot\normal_{FT}
			+\dfrac{s_{T'}s_{T}}{s_{T'}-s_{T}}\dfrac{1}{a^2}\left(p_{T'}  -p_{T}\right).
\end{align*}
The two external acoustic wave speeds are
\begin{equation*}
	s_{T'}	=	\dfrac{1}{2}
				\left(
				\uVec_{T'}\cdot\normal_{TF}
				+\sqrt{\left(\uVec_{T'}\cdot\normal_{TF}\right)^{2}  +4a^{2}}
				\right)
	\qquad
	\text{and}
	\qquad
	s_{T}	=	\dfrac{1}{2}
				\left(
				\uVec_{T}\cdot\normal_{TF}
				-\sqrt{\left(\uVec_{T}\cdot\normal_{TF}\right)^{2}  +4a^{2}}
				\right),
\end{equation*}
regardless of the compressive or expansive nature of the acoustic waves, see Appendix B4 in \cite{Bassi.Massa.ea:18} for details.
Defining the intermediate state as
\begin{equation*}
	\uVec_{F}	=	\dfrac{s_{T'}\uVec_{T}  -s_{T}\uVec_{T'}}{s_{T'}-s_{T}}
	\qquad
	\text{and}
	\qquad
	p_{F}	=	\dfrac{s_{T'}p_{T}  -s_{T}p_{T'}}{s_{T'}-s_{T}},
\end{equation*}
and performing some trivial algebraic manipulations, the flux${_\star}$ can be written as follows
\begin{align*}
	\left(\uVec\otimes\uVec  +p\I \right)_{\star} \cdot\normal_{FT}
		&=	\begin{aligned}[t]
			&\left(\uVec_{F}\otimes\uVec_{T}  +p_{F}\I\right) \cdot\normal_{FT}
			+s_{T'}\left(\uVec_{T}  -\uVec_{F}\right)
			-\dfrac{s_{T}}{s_{T'}-s_{T}}
				\uVec_{T'}
				\left[\left(\uVec_{T'}  -\uVec_{T}\right) \cdot \normal_{TF}\right],
			\end{aligned}\\
	\left(\uVec\right)_{\star} \cdot\normal_{FT}
		&=	\uVec_{F} \cdot\normal_{FT}
			+s_{T'}\dfrac{1}{a^2}\left(p_{T}  -p_{F}\right).
\end{align*}
In order to decouple neighboring elements, the continuity of the normal component of the velocity is assumed by hypothesis, 
namely $\uVec_{T}\cdot\normal_{TF}=\uVec_{T'}\cdot\normal_{TF}=\uVec_{F}\cdot\normal_{TF}$.
Thus, the numerical fluxes in Equations \eqref{eq:HLL_numerical_flux_qdm}-\eqref{eq:HLL_numerical_flux_cnt} 
as well as the interface conditions in \eqref{eq:face_variables} and the wave speeds in \eqref{eq:waveSpeed} are directly obtained. 
We remark that, according to the HLL machinery, see \cite{Toro:2009} for details, unique definition of numerical fluxes follows from the observation that $s_{T'}>0$ and $s_{T}<0$.

\section{Numerical results}
\label{sec:numResults}

\subsection{Kovasznay flow}
\label{sec:Kova}
Steady flow behind a grid made of equally spaced parallel rods
is described by the following exact solution of the incompressible Navier--Stokes equations (see \cite{Kovasznay:48})
\begin{equation}\label{eq:kovasznay_solution}
    \begin{aligned}
       \uVec 
            &= \left[1-\eulern^{\kappa x}\cos{\rbrackets{2\pi y}}\right] \iVec +
             \dfrac{\kappa }{2\pi}\eulern^{\kappa x}\sin{\rbrackets{2\pi y}} \jVec,
            \\
        p
	    &=  p_{0}-\dfrac{1}{2}\eulern^{2\kappa x},
    \end{aligned}
\end{equation}
where $p_{0}\in\Real$ is an arbitrary constant and the parameter $\kappa$ depends on the Reynolds number  
\begin{equation*}
    \kappa  =   \dfrac{\Reynolds}{2}    -\sqrt{\dfrac{\Reynolds^{2}}{4} +4\pi^{2}}.
\end{equation*}
The Kovasznay flow problem is defined on the 
bi-unit square computational domain $\Omega = \rbrackets{-0.5,1.5}{\times}\rbrackets{0,2}$
and solved by means of the \hhohdiv and \hhohll schemes. 
Boundary conditions are derived from the Kovasznay analytical solution \eqref{eq:kovasznay_solution} at $\Reynolds = \frac{1}{\nu} = 40$.
In particular, a Neumann boundary condition is imposed at the outflow boundary (right side of the square domain) while Dirichlet boundary conditions are imposed on the remaining boundaries. 

In order to numerically validate the proposed HHO schemes,
we consider $h$-refined mesh sequences composed of regular triangular and quadrilateral elements, 
doubling the number of mesh elements in each Cartesian direction at each refinement step, and several polynomial degrees $k=\{1,2,3,4\}$.
Errors in $L^2$-norm and $h$-convergence rates are tabulated for the approximated 
pressure, divergence of the velocity, velocity, and velocity gradients fields. 
Results over triangular elements meshes are reported in Tables~\ref{tab:kovasznay_Hdiv_tri} and \ref{tab:kovasznay_HLL_tri} for \hhohdiv and \hhohll schemes, respectively.
Results over quadrilateral elements meshes are reported in Tables~\ref{tab:kovasznay_Hdiv_quad} and \ref{tab:kovasznay_HLL_quad} for \hhohdiv and \hhohll schemes, respectively.
As a general observation it can be noticed that, although the level of accuracy provided by the two schemes proposed is comparable, the \hhohdiv scheme takes the lead over triangular elements meshes while
\hhohll is to be preferred over quadrilateral meshes.

Since the flow is diffusion dominated at this Reynolds number, the \hhohll scheme shows an asymptotic 
convergence rates of $k+2$ for the velocity and $k+1$ for pressure, velocity divergence and velocity gradients.
The convergence rates obtained over quadrilateral elements meshes are slightly better than those on triangular elements meshes, at all polynomial degrees. 
Notice that, at the same mesh sequence refinement step, quadrilateral meshes have twice the number of elements as compared to triangular meshes.

\hhohdiv scheme shows asymptotic convergence rates of order $k+1$ for the pressure unknown on both triangular and quadrilateral meshes.  
The incompressibility constraint is exactly satisfied, note that the velocity divergence errors settle around the machine precision.  
Convergence rates for the velocity and the velocity gradients differs when applying \hhohdiv on triangular and quadrilateral elements meshes. 
On triangular meshes we observe the optimal rates of $k+2$ and $k+1$ for the velocity and the velocity gradients, respectively, while 
on quadrilateral meshes the convergence tops at $k+1$ and $k$, respectively.
Accordingly, a full order of convergence is lost.
We refer to the work of Kirk \ea\ \cite{Kirk2019} for additional clues regarding improved convergence rates on simplexes.

Eventually, as remarked in Section \ref{subsec:hdivConf}, we verify that optimal convergence rates on quadrilateral meshes can 
be recovered by lowering the degree of polynomials employed to approximate the numerical trace of the pressure on mesh faces.
Results for the \hhodivfree scheme on the quadrilateral mesh sequence are reported in Table~\ref{tab:kovasznay_Hnodiv_quad}.
Since the \hhodivfree scheme is not H-div conforming and has demonstrated to be less robust than \hhohll in convection-dominated flow regimes, 
it will not be further investigated in subsequent test cases.

\begin{table}[!t]
\centering
\caption{Kovasznay flow. Errors and $h$-convergence rates for the \hhohdiv scheme on the triangular elements mesh sequence.}
\label{tab:kovasznay_Hdiv_tri}
\begin{tabular}{cccccccccc}
\hline
&
$\card{\Th}$  &
$\left\Vert \nabla \uVec_{T}  -\nabla \uVec \right\Vert_{L^2}$ &
order &
$\left\Vert \uVec_{T}  -\uVec \right\Vert_{L^2}$ &
order &
$\left\Vert p_{T}  -p \right\Vert_{L^2}$ &
order &
$\left\Vert  \nabla \cdot \uVec_{T} \right\Vert_{L^2}$ &
order \\
\hline
\multirow{5}{*}{\rotatebox{90}{$k=1$}} &       32      &       2.42e+00        &       --      &       1.08e-01        &       --      &       5.68e-02        &       --      &       7.32e-16        &       --       \\
        &       128     &       8.90e-01        &       1.44    &       2.58e-02        &       2.07    &       1.82e-02        &       1.64    &       5.74e-16        &       --       \\
        &       512     &       2.99e-01        &       1.57    &       3.63e-03        &       2.83    &       4.99e-03        &       1.87    &       4.74e-16        &       --       \\
        &       2048    &       8.51e-02        &       1.81    &       4.96e-04        &       2.87    &       1.18e-03        &       2.08    &       5.19e-16        &       --       \\
        &       8192    &       2.28e-02        &       1.90    &       6.65e-05        &       2.90    &       2.87e-04        &       2.03    &       5.09e-16        &       --		\\
\hline
\multirow{5}{*}{\rotatebox{90}{$k=2$}}        &       32      &        9.24e-01   & --         & 3.78e-02   & --         & 1.52e-02   & --         & 6.60e-16   & --   \\
        &       128     &        1.60e-01   & 2.52       & 3.05e-03   & 3.63       & 1.46e-03   & 3.38       & 6.91e-16   & --   \\
        &       512     &        2.33e-02   & 2.78       & 1.72e-04   & 4.15       & 1.97e-04   & 2.89       & 5.52e-16   & --   \\
        &       2048    &        3.80e-03   & 2.61       & 1.46e-05   & 3.55       & 2.73e-05   & 2.85       & 5.68e-16   & --   \\
        &       8192    &        4.75e-04   & 3.00       & 8.96e-07   & 4.03       & 3.08e-06   & 3.15       & 5.48e-16   & --   \\
\hline
\multirow{5}{*}{\rotatebox{90}{$k=3$}}        &       32      &       1.95e-01   & --         & 4.85e-03   & --         & 2.13e-03   & --         & 6.73e-16    & --   \\
        &       128     &       2.10e-02   & 3.21       & 2.69e-04   & 4.17       & 1.75e-04   & 3.60       & 6.88e-16    & --   \\
        &       512     &       1.83e-03   & 3.52       & 1.02e-05   & 4.73       & 1.23e-05   & 3.84       & 5.45e-16    & --   \\
        &       2048    &       1.24e-04   & 3.89       & 3.43e-07   & 4.89       & 6.67e-07   & 4.20       & 5.99e-16    & --   \\
        &       8192    &       8.14e-06   & 3.92       & 1.12e-08   & 4.93       & 4.04e-08   & 4.04       & 5.99e-16    & --   \\
\hline
\multirow{5}{*}{\rotatebox{90}{$k=4$}}        &       32      &       4.44e-02   & --         & 9.20e-04   & --         & 3.63e-04   & --         & 7.73e-16    & --   \\
        &       128     &       2.02e-03   & 4.46       & 2.51e-05   & 5.19       & 1.10e-05   & 5.04       & 7.01e-16    & --   \\
        &       512     &       7.67e-05   & 4.72       & 3.19e-07   & 6.30       & 3.65e-07   & 4.92       & 6.03e-16    & --   \\
        &       2048    &       3.32e-06   & 4.53       & 7.31e-09   & 5.45       & 1.35e-08   & 4.76       & 6.33e-16    & --   \\
        &       8192    &       9.90e-08   & 5.07       & 1.05e-10   & 6.12       & 3.84e-10   & 5.13       & 6.30e-16    & --
\end{tabular}
\end{table}
\begin{table}[!t]
\centering
\caption{Kovasznay flow. Errors and $h$-convergence rates for the \hhohll scheme on the triangular element mesh sequence.}
\label{tab:kovasznay_HLL_tri}
\begin{tabular}{cccccccccc}
\hline
&
$\card{\Th}$  &
$\left\Vert \nabla \uVec_{T}  -\nabla \uVec \right\Vert_{L^2}$ &
order &
$\left\Vert \uVec_{T}  -\uVec \right\Vert_{L^2}$ &
order &
$\left\Vert p_{T}  -p \right\Vert_{L^2}$ &
order &
$\left\Vert  \nabla \cdot \uVec_{T} \right\Vert_{L^2}$ &
order \\
\hline
\multirow{5}{*}{\rotatebox{90}{$k=1$}} &       32      &       3.22e+00   & --         & 3.30e-01   & --         & 2.38e-01   & --         & 4.80e-01    & --   \\
        &       128     &       1.46e+00   & 1.22       & 7.54e-02   & 2.28       & 4.70e-02   & 2.50       & 3.62e-01    & 0.41   \\
        &       512     &       5.23e-01   & 1.45       & 1.30e-02   & 2.49       & 7.56e-03   & 2.59       & 1.42e-01    & 1.35   \\
        &       2048    &       1.66e-01   & 1.71       & 1.89e-03   & 2.88       & 1.38e-03   & 2.55       & 5.36e-02    & 1.41   \\
        &       8192    &       4.76e-02   & 1.79       & 2.66e-04   & 2.81       & 2.65e-04   & 2.36       & 1.75e-02    & 1.61   \\
\hline
\multirow{5}{*}{\rotatebox{90}{$k=2$}} &       32      &       1.37e+00   & --         & 8.12e-02   & --         & 6.24e-02   & --         & 2.65e-01    & --   \\
        &       128     &       2.10e-01   & 2.89       & 7.22e-03   & 3.73       & 4.67e-03   & 3.99       & 4.87e-02    & 2.44   \\
        &       512     &       3.20e-02   & 2.67       & 4.95e-04   & 3.80       & 3.83e-04   & 3.54       & 8.45e-03    & 2.53   \\
        &       2048    &       5.29e-03   & 2.69       & 4.32e-05   & 3.64       & 5.34e-05   & 2.95       & 1.47e-03    & 2.52   \\
        &       8192    &       7.22e-04   & 2.85       & 2.97e-06   & 3.84       & 6.94e-06   & 2.92       & 2.21e-04    & 2.73   \\
\hline
\multirow{5}{*}{\rotatebox{90}{$k=3$}} &       32      &       2.22e-01   & --         & 1.06e-02   & --         & 3.99e-03   & --         & 3.97e-02    & --   \\
        &       128     &       2.23e-02   & 3.54       & 5.15e-04   & 4.65       & 2.70e-04   & 4.15       & 4.93e-03    & 3.01   \\
        &       512     &       1.96e-03   & 3.45       & 2.41e-05   & 4.34       & 1.95e-05   & 3.72       & 4.63e-04    & 3.41   \\
        &       2048    &       1.44e-04   & 3.90       & 8.77e-07   & 4.95       & 1.41e-06   & 3.93       & 4.01e-05    & 3.53   \\
        &       8192    &       9.91e-06   & 3.84       & 3.08e-08   & 4.80       & 9.03e-08   & 3.94       & 2.95e-06    & 3.76   \\
\hline
\multirow{5}{*}{\rotatebox{90}{$k=4$}} &       32      &       4.44e-02   & --         & 1.63e-03   & --         & 1.10e-03   & --         & 1.14e-02    & --   \\
        &       128     &       1.78e-03   & 4.95       & 3.24e-05   & 6.04       & 2.27e-05   & 5.98       & 4.42e-04    & 4.69   \\
        &       512     &       6.99e-05   & 4.59       & 6.46e-07   & 5.54       & 8.01e-07   & 4.74       & 1.93e-05    & 4.52   \\
        &       2048    &       3.02e-06   & 4.70       & 1.47e-08   & 5.66       & 3.12e-08   & 4.85       & 8.39e-07    & 4.52   \\
        &       8192    &       1.05e-07   & 4.81       & 2.62e-10   & 5.77       & 1.09e-09   & 4.82       & 3.07e-08    & 4.77
\end{tabular}
\end{table}
\begin{table}[!t]
\centering
\caption{Kovasznay flow. Errors and $h$-convergence rates for the \hhohdiv scheme on the quadrilateral elements mesh sequence.}
\label{tab:kovasznay_Hdiv_quad}
\begin{tabular}{cccccccccc}
\hline
&
$\card{\Th}$  &
$\left\Vert \nabla \uVec_{T}  -\nabla \uVec \right\Vert_{L^2}$ &
order &
$\left\Vert \uVec_{T}  -\uVec \right\Vert_{L^2}$ &
order &
$\left\Vert p_{T}  -p \right\Vert_{L^2}$ &
order &
$\left\Vert  \nabla \cdot \uVec_{T} \right\Vert_{L^2}$ &
order \\
\hline
\multirow{5}{*}{\rotatebox{90}{$k=1$}}  &       64      &       4.01e+00   & --         & 1.40e-01   & --         & 2.93e-02   & --         & 3.91e-16    & --   \\
        &       256     &       1.99e+00   & 1.01       & 3.53e-02   & 1.99       & 7.96e-03   & 1.88       & 3.14e-16    & --   \\
        &       1024    &       9.90e-01   & 1.00       & 8.97e-03   & 1.98       & 2.27e-03   & 1.81       & 3.18e-16    & --   \\
        &       4096    &       4.94e-01   & 1.00       & 2.26e-03   & 1.99       & 5.98e-04   & 1.93       & 2.79e-16    & --   \\
        &       16384   &       2.47e-01   & 1.00       & 5.67e-04   & 1.99       & 1.52e-04   & 1.97       & 3.04e-16    & --   \\
\hline
\multirow{5}{*}{\rotatebox{90}{$k=2$}}  &       64      &       8.88e-01   & --         & 1.96e-02   & --         & 4.01e-03   & --         & 3.12e-16    & --   \\
        &       256     &       2.15e-01   & 2.05       & 2.33e-03   & 3.07       & 3.31e-04   & 3.60       & 3.67e-16    & --   \\
        &       1024    &       5.32e-02   & 2.01       & 2.84e-04   & 3.04       & 2.36e-05   & 3.81       & 3.42e-16    & --   \\
        &       4096    &       1.33e-02   & 2.00       & 3.52e-05   & 3.01       & 1.78e-06   & 3.73       & 3.34e-16    & --   \\
        &       16384   &       3.31e-03   & 2.00       & 4.39e-06   & 3.00       & 1.59e-07   & 3.49       & 3.75e-16    & --   \\
\hline
\multirow{5}{*}{\rotatebox{90}{$k=3$}}  &       64      &       1.26e-01   & --         & 1.82e-03   & --         & 2.67e-04   & --         & 3.32e-16    & --   \\
        &       256     &       1.53e-02   & 3.04       & 1.15e-04   & 3.98       & 2.29e-05   & 3.54       & 3.63e-16    & --   \\
        &       1024    &       1.90e-03   & 3.01       & 7.22e-06   & 3.99       & 1.57e-06   & 3.87       & 3.60e-16    & --   \\
        &       4096    &       2.36e-04   & 3.00       & 4.52e-07   & 4.00       & 1.00e-07   & 3.97       & 3.39e-16    & --   \\
        &       16384   &       2.95e-05   & 3.00       & 2.83e-08   & 4.00       & 6.31e-09   & 3.99       & 3.70e-16    & --   \\
\hline
\multirow{5}{*}{\rotatebox{90}{$k=4$}}  &       64      &       1.35e-02   & --         & 1.46e-04   & --         & 2.12e-05   & --         & 3.29e-16    & --   \\
        &       256     &       8.16e-04   & 4.05       & 4.32e-06   & 5.08       & 4.33e-07   & 5.61       & 3.67e-16    & --   \\
        &       1024    &       5.04e-05   & 4.02       & 1.32e-07   & 5.03       & 8.65e-09   & 5.64       & 4.06e-16    & --   \\
        &       4096    &       3.14e-06   & 4.00       & 4.11e-09   & 5.01       & 2.00e-10   & 5.44       & 3.94e-16    & --   \\
        &       16384   &       1.96e-07   & 4.00       & 1.28e-10   & 5.00       & 5.40e-12   & 5.21       & 4.30e-16    & --
\end{tabular}
\end{table}
\begin{table}[!t]
\centering
\caption{Kovasznay flow. Errors and $h$-convergence rates for the \hhohll scheme on the quadrilateral elements mesh sequence.}
\label{tab:kovasznay_HLL_quad}
\begin{tabular}{cccccccccc}
\hline
&
$\card{\Th}$  &
$\left\Vert \nabla \uVec_{T}  -\nabla \uVec \right\Vert_{L^2}$ &
order &
$\left\Vert \uVec_{T}  -\uVec \right\Vert_{L^2}$ &
order &
$\left\Vert p_{T}  -p \right\Vert_{L^2}$ &
order &
$\left\Vert  \nabla \cdot \uVec_{T} \right\Vert_{L^2}$ &
order \\
\hline
\multirow{5}{*}{\rotatebox{90}{$k=1$}}  &       64      &       2.11e+00   & --         & 1.08e-01   & --         & 7.82e-02   & --         & 1.38e-01    & --   \\
        &       256     &       7.63e-01   & 1.53       & 1.67e-02   & 2.82       & 1.20e-02   & 2.82       & 8.86e-02    & 0.64   \\
        &       1024    &       2.45e-01   & 1.67       & 2.69e-03   & 2.69       & 1.97e-03   & 2.67       & 3.44e-02    & 1.36   \\
        &       4096    &       7.19e-02   & 1.79       & 4.12e-04   & 2.74       & 3.11e-04   & 2.69       & 1.07e-02    & 1.68   \\
        &       16384   &       1.97e-02   & 1.88       & 5.82e-05   & 2.84       & 5.08e-05   & 2.63       & 2.99e-03    & 1.84   \\
\hline
\multirow{5}{*}{\rotatebox{90}{$k=2$}}  &       64      &       3.42e-01   & --         & 1.39e-02   & --         & 6.62e-03   & --         & 3.70e-02    & --   \\
        &       256     &       5.94e-02   & 2.64       & 1.28e-03   & 3.59       & 7.09e-04   & 3.36       & 7.78e-03    & 2.25   \\
        &       1024    &       9.17e-03   & 2.75       & 9.73e-05   & 3.80       & 5.73e-05   & 3.71       & 1.39e-03    & 2.48   \\
        &       4096    &       1.29e-03   & 2.86       & 6.77e-06   & 3.89       & 5.61e-06   & 3.39       & 2.13e-04    & 2.71   \\
        &       16384   &       1.73e-04   & 2.92       & 4.48e-07   & 3.94       & 6.49e-07   & 3.13       & 2.97e-05    & 2.84   \\
\hline
\multirow{5}{*}{\rotatebox{90}{$k=3$}}  &       64      &       4.17e-02   & --         & 1.23e-03   & --         & 5.86e-04   & --         & 7.60e-03    & --   \\
        &       256     &       3.25e-03   & 3.84       & 5.05e-05   & 4.80       & 2.54e-05   & 4.72       & 7.69e-04    & 3.30   \\
        &       1024    &       2.35e-04   & 3.87       & 1.92e-06   & 4.82       & 1.28e-06   & 4.41       & 6.43e-05    & 3.58   \\
        &       4096    &       1.61e-05   & 3.91       & 6.73e-08   & 4.89       & 5.72e-08   & 4.53       & 4.85e-06    & 3.73   \\
        &       16384   &       1.06e-06   & 3.94       & 2.24e-09   & 4.94       & 2.66e-09   & 4.45       & 3.38e-07    & 3.84   \\
\hline
\multirow{5}{*}{\rotatebox{90}{$k=4$}}  &       64      &       3.73e-03   & --         & 8.16e-05   & --         & 4.07e-05   & --         & 5.57e-04    & --   \\
        &       256     &       1.49e-04   & 4.84       & 1.59e-06   & 5.93       & 1.12e-06   & 5.41       & 2.01e-05    & 4.79   \\
        &       1024    &       5.38e-06   & 4.90       & 2.80e-08   & 5.95       & 3.29e-08   & 5.20       & 7.31e-07    & 4.78   \\
        &       4096    &       1.81e-07   & 4.94       & 4.67e-10   & 5.97       & 9.85e-10   & 5.12       & 2.57e-08    & 4.83   \\
        &       16384   &       5.91e-09   & 4.97       & 7.57e-12   & 5.98       & 2.97e-11   & 5.08       & 8.95e-10    & 4.84
\end{tabular}
\end{table}

\begin{table}[!t]
\centering
\caption{Kovasznay flow. Errors and $h$-convergence rates for the \hhodivfree scheme on the quadrilateral elements mesh sequence.}
\label{tab:kovasznay_Hnodiv_quad}
\begin{tabular}{cccccccccc}
\hline
&
$\card{\Th}$  &
$\left\Vert \nabla \uVec_{T}  -\nabla \uVec \right\Vert_{L^2}$ &
order &
$\left\Vert \uVec_{T}  -\uVec \right\Vert_{L^2}$ &
order &
$\left\Vert p_{T}  -p \right\Vert_{L^2}$ &
order &
$\left\Vert  \nabla \cdot \uVec_{T} \right\Vert_{L^2}$ &
order \\
\hline
\multirow{5}{*}{\rotatebox{90}{$k=1$}}  
        &       64      & 1.26e+00   & --         & 5.65e-02   & --         & 2.11e-02   & --         & 1.73e-16       & --   \\
        &       256     & 3.07e-01   & 2.04       & 6.25e-03   & 3.18       & 2.33e-03   & 3.18       & 2.23e-16       & --   \\
        &       1024    & 7.63e-02   & 2.01       & 7.80e-04   & 3.00       & 5.17e-04   & 2.17       & 2.92e-16       & --   \\
        &       4096    & 1.91e-02   & 2.00       & 9.83e-05   & 2.99       & 1.28e-04   & 2.02       & 3.31e-16       & --   \\
        &       16384   & 4.79e-03   & 2.00       & 1.24e-05   & 2.99       & 3.19e-05   & 2.00       & 3.17e-16       & --   \\
\hline
\multirow{5}{*}{\rotatebox{90}{$k=2$}}  
        &       64      &     2.63e-01   & --         & 6.64e-03   & --         & 2.03e-03   & --         & 2.02e-16   & --   \\
        &       256     &     3.46e-02   & 2.93       & 4.53e-04   & 3.87       & 1.78e-04   & 3.51       & 2.72e-16   & --   \\
        &       1024    &     4.42e-03   & 2.97       & 2.91e-05   & 3.96       & 1.42e-05   & 3.65       & 2.95e-16   & --   \\
        &       4096    &     5.58e-04   & 2.98       & 1.84e-06   & 3.98       & 1.26e-06   & 3.49       & 3.69e-16   & --   \\
        &       16384   &     7.02e-05   & 2.99       & 1.16e-07   & 3.99       & 1.32e-07   & 3.26       & 3.44e-16   & --   \\
\hline
\multirow{5}{*}{\rotatebox{90}{$k=3$}}  
        &       64      &       3.93e-02   & --         & 6.28e-04   & --         & 1.58e-04   & --         & 1.97e-16      & --   \\
        &       256     &       2.64e-03   & 3.90       & 2.07e-05   & 4.93       & 4.46e-06   & 5.15       & 2.92e-16      & --   \\
        &       1024    &       1.70e-04   & 3.96       & 6.67e-07   & 4.95       & 1.79e-07   & 4.64       & 3.40e-16      & --   \\
        &       4096    &       1.08e-05   & 3.98       & 2.11e-08   & 4.98       & 7.91e-09   & 4.50       & 3.86e-16      & --   \\
        &       16384   &       6.77e-07   & 3.99       & 6.65e-10   & 4.99       & 3.64e-10   & 4.44       & 3.71e-16      & --   \\
\hline
\multirow{5}{*}{\rotatebox{90}{$k=4$}}  
        &       64      &       4.43e-03   & --         & 5.07e-05   & --         & 1.05e-05   & --         & 2.39e-16   & --   \\
        &       256     &       1.49e-04   & 4.90       & 8.67e-07   & 5.87       & 2.61e-07   & 5.33       & 3.05e-16   & --   \\
        &       1024    &       4.76e-06   & 4.96       & 1.40e-08   & 5.96       & 6.45e-09   & 5.34       & 3.49e-16   & --   \\
        &       4096    &       1.50e-07   & 4.98       & 2.21e-10   & 5.98       & 1.76e-10   & 5.20       & 4.14e-16   & --   \\
        &       16384   &       4.73e-09   & 4.99       & 3.48e-12   & 5.99       & 5.13e-12   & 5.10       & 3.98e-16   & --
\end{tabular}
\end{table}
\subsection{LLMS pressure gradient test case}
\label{sec:LLMS}
In order to investigate pressure-robustness, \ie velocity solution accuracy in the presence of strong pressure gradients \cite{Lederer.Linke.ea:17,Rhebergen.Wells:18},
we consider the test case proposed by Lederer, Linke, Merdon and Sch\"{o}berl (LLMS) \cite{Lederer.Linke.ea:17},
hereinafter referred to as LLMS pressure gradient test case.
We impose boundary conditions and forcing term $\boldsymbol{f}$ according to 
the following analytical velocity and pressure fields
\begin{equation*}
	\uVec	=	\nabla\times\zeta
    \qquad
    \text{and}
    \qquad
    p	=	p_{0}  +x^{7}  +y^{7},
\end{equation*}
with arbitrary $p_{0}\in\Real$ and $\zeta = x^{2}\left(x-1\right)^{2}y^{2}\left(y-1\right)^{2}$. 
In case of \hhohdiv we impose Dirichlet boundary conditions on all but one side of the unit square domain $\Omega = \left(0,1\right)\times\left(0,1\right)$, where a Neumann boundary is set.
In order to ensure stability in the inviscid limit, in case of \hhohll we impose both Dirichlet and Neumann boundary conditions on all sides, leading to the so called \emph{given} boundary conditions.
We remark that, according to given BCs, HLL Riemann problems on boundary faces utilize both the analytical pressure and velocity solutions for the definition of external states.

Figure~\ref{fig:analytical_pressure-robustness} shows the horizontal velocity component ($\uVec \cdot \iVec$) and pressure behavior over the computational domain.
Notice the presence of a strong pressure gradient near the top-right corner of the square.
\begin{figure}[!t]
    \centering
	\includegraphics[width=0.4\textwidth]{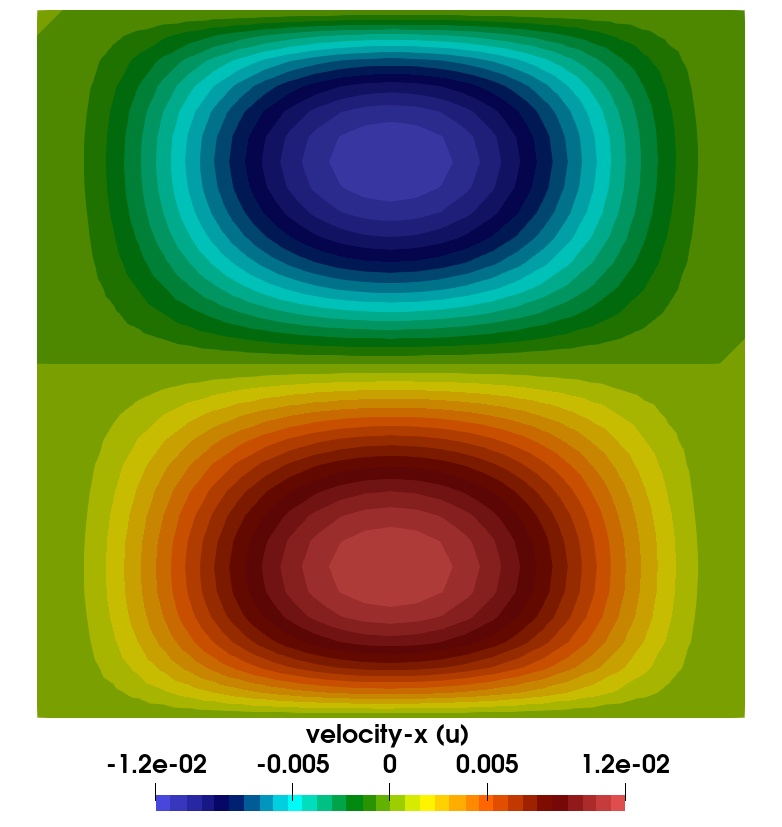}
	\includegraphics[width=0.4\textwidth]{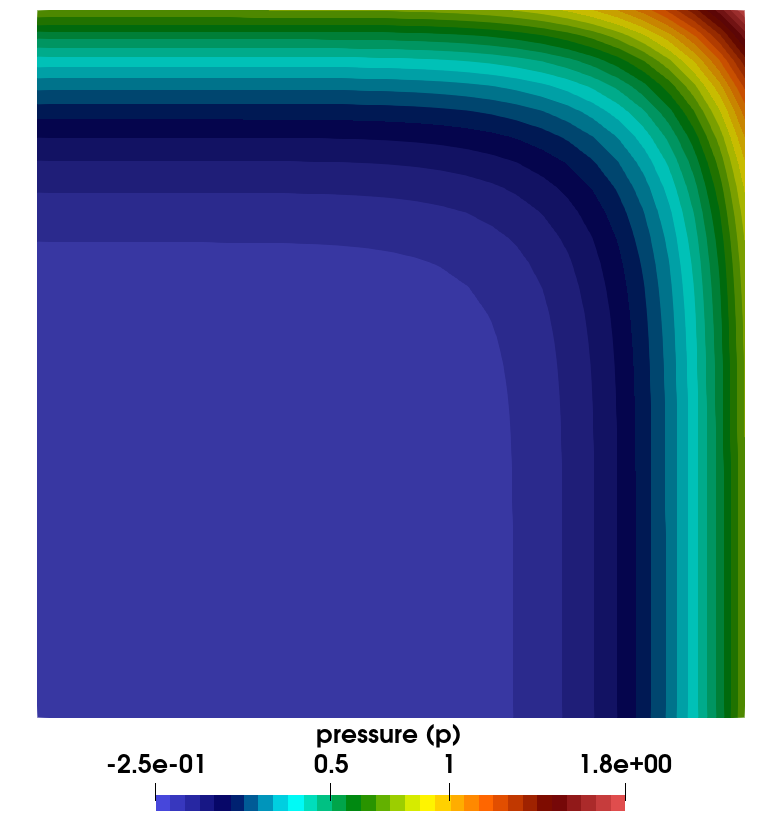}
    \caption{LLMS pressure gradient. Horizontal velocity component ($u$) and pressure ($p$) fields ($p_{0}=-0.25$).}
	\label{fig:analytical_pressure-robustness}
\end{figure}

In order to numerically validate convergence rates in the diffusion dominated and convection dominated flow regimes we consider $\nu=1$ and $\nu=10^{-4}$, respectively.
The computational domain is discretized by means of increasingly fine regular triangular elements meshes, doubling the number of elements in each Cartesian direction at each refinement step, and we focus on $k=3$ HHO formulations.
Tables \ref{tab:pressure_robustness_Hdiv_ni_1}-\ref{tab:pressure_robustness_HLL_ni_10-4} summarize the convergence analysis 
considering errors in $L^2$-norm for the velocity, the velocity gradients and the pressure fields. 
In the diffusion dominated regime both \hhohdiv and \hhohll deliver the same convergence rates, \ie order $k+2$ for the velocity and $k+1$ for pressure and velocity gradients.
In the convection dominated regime, only \hhohdiv maintains the aforementioned convergence rates, while the \hhohll scheme loses an order on velocity and velocity gradients fields.
\begin{table}[!t]
\centering
\caption{LLMS pressure gradient. Errors and $h$-convergence rates for $k=3$ \hhohdiv with $\nu=1$.}
\label{tab:pressure_robustness_Hdiv_ni_1}
\begin{tabular}{ccccccc}
\hline
$\card{\Th}$  &
$\left\Vert \nabla \uVec_{T}  -\nabla \uVec \right\Vert_{L^2}$ &
order &
$\left\Vert \uVec_{T}  -\uVec \right\Vert_{L^2}$ &
order &
$\left\Vert p_{T}  -p \right\Vert_{L^2}$ &
order \\
\hline
32      &   7.22e-04   	&	--      &	1.02e-05    &   --      &   3.07e-04    &   --      \\ 
128     &   5.34e-05   	& 	3.76    &	4.05e-07    &   4.65    &   1.93e-05    &   3.99    \\ 
512     &   3.60e-06   	& 	3.89    &	1.43e-08    &   4.83    &   1.18e-06    &   4.04    \\ 
2048    &   2.34e-07   	& 	3.95    &	4.74e-10    &   4.91    &   7.20e-08    &   4.03    \\ 
8192    &   1.49e-08	& 	3.97 	&	1.53E-11    &   4.96    &   4.43e-09    &   4.02     
\end{tabular}
\end{table}
\begin{table}[!t]
\centering
\caption{LLMS pressure gradient. Errors and $h$-convergence rates for $k=3$ \hhohll with $\nu=1$.}
\begin{tabular}{ccccccc}
\hline
$\card{\Th}$  &
$\left\Vert \nabla \uVec_{T}  -\nabla \uVec \right\Vert_{L^2}$ &
order &
$\left\Vert \uVec_{T}  -\uVec \right\Vert_{L^2}$ &
order &
$\left\Vert p_{T}  -p \right\Vert_{L^2}$ &
order \\
\hline
32      &   4.59e-04   	& 	--      &	1.14e-05    &   --      &   4.50e-04    &   --      \\ 
128     &   2.85e-05  	& 	4.01    &	3.63e-07    &   4.97    &   2.94e-05    &   3.93    \\ 
512     &   1.72e-06   	& 	4.05    &	1.14e-08    &   4.99    &   1.86e-06    &   3.99    \\ 
2048    &   1.04e-07   	& 	4.04    &	3.59e-10    &   4.99    &   1.16e-07    &   4.00    \\ 
8192    &   6.41e-09   	& 	4.03    &	1.13e-11    &   4.99    &   7.23e-09    &   4.00     
\end{tabular}
\end{table}
\begin{table}[!t]
\centering
\caption{LLMS pressure gradient. Errors and $h$-convergence rates for $k=3$ \hhohdiv with $\nu=10^{-4}$.}
\begin{tabular}{ccccccc}
\hline
$\card{\Th}$  &
$\left\Vert \nabla \uVec_{T}  -\nabla \uVec \right\Vert_{L^2}$ &
order &
$\left\Vert \uVec_{T}  -\uVec \right\Vert_{L^2}$ &
order &
$\left\Vert p_{T}  -p \right\Vert_{L^2}$ &
order \\
\hline
32      &   1.59e-03   & --         & 1.86e-05   & --         & 2.58e-04 & --		\\
128     &   8.18e-05   & 4.28       & 5.51e-07   & 5.07       & 1.67e-05   & 3.95	\\
512     &   6.62e-06   & 3.63       & 2.28e-08   & 4.60       & 1.05e-06   & 3.99	\\
2048    &   4.72e-07   & 3.81       & 8.30e-10   & 4.78       & 6.60e-08   & 4.00 	\\
8192    &   3.09e-08   & 3.93       & 2.75e-11   & 4.92       & 4.13e-09   & 4.00
\end{tabular}
\end{table}   
\begin{table}[!t]
\centering
\caption{LLMS pressure gradient. Errors and $h$-convergence rates for $k=3$ \hhohll with $\nu=10^{-4}$.}
\label{tab:pressure_robustness_HLL_ni_10-4}
\begin{tabular}{ccccccc}
\hline
$\card{\Th}$  &
$\left\Vert \nabla \uVec_{T}  -\nabla \uVec \right\Vert_{L^2}$ &
order &
$\left\Vert \uVec_{T}  -\uVec \right\Vert_{L^2}$ &
order &
$\left\Vert p_{T}  -p \right\Vert_{L^2}$ &
order \\
\hline
32      &   2.15e-02   & 	--         &	5.77e-04    &   --      &   2.65e-04    &   --      \\ 
128     &   2.22e-03   & 	3.28       &	2.45e-05    &   4.56    &   1.72e-05    &   3.95    \\ 
512     &   2.50e-04   & 	3.15       &	1.24e-06    &   4.31    &   1.09e-06    &   3.98    \\ 
2048    &   2.88e-05   & 	3.12       &	6.85e-08    &   4.18    &   6.82e-08    &   3.99    \\ 
8192    &   3.24e-06   & 	3.15       &	3.88e-09    &   4.14    &   4.27e-09    &   4.00  
\end{tabular}
\end{table}

Figure~\ref{fig:pressure_robustness_nu_analysis} depicts the pressure and velocity errors in $L^2$-norm over the 2k triangular elements mesh 
for different viscosities values, \ie $\nu = \{1,10^{-1},10^{-2},10^{-3},10^{-4}\}$.
The fact that velocity errors are almost unaltered while varying the viscosity confirms that the \hhohdiv formulation is pressure-robust. 
Note that the velocity errors provided by the \hhohll scheme increase by two orders of magnitude while decreasing the viscosity by four orders of magnitude.

We next investigate robustness in the inviscid limit considering $k=3$ and $k=4$ HHO formulations over triangular and quadrilateral meshes 
and varying the viscosity in one order of magnitude steps from $1$ to $10^{-14}$.
The error analysis reported in Figures~\ref{fig:pressure_robustness_inviscid_p3} and \ref{fig:pressure_robustness_inviscid_p4} for different mesh densities, see figures captions for details,
shows that pressure-robustness does not guarantee robustness in the inviscid limit. 
Indeed, only the \hhohll formulation is stable in the limit of vanishing viscosity.
The velocity errors provided by the \hhohdiv formulation tend to dramatically increase for $\nu < 10^{-5}$. 
Moreover, due to convergence failure of Newton's method globalization strategy \cite{Kelley.Keyes:98}, the numerical solutions are not available for $\nu < 10^{-7}$. 
We verified that even using the $L^2$-projection of the exact solution as initial guess the simulation blows up.
Interestingly, the \hhohll formulation shows to be resilient in the inviscid limit: notice that the velocity error reaches a plateau while decreasing the 
viscosity below $10^{-10}$ and also the pressure error is well behaved.  
In case of $k=4$ \hhohll formulations over the 8k triangular and the 4k quadrilateral mesh, the velocity errors increase by three orders of magnitude
moving from the diffusion dominated regime to the inviscid limit, while the pressure error stays almost constant around $10^{-11}$. 

\begin{figure} [!t]
\centering
\begin{tikzpicture}[scale=0.65]
\begin{axis}    [
                unbounded coords=jump,
                ymode=log,xmode=log,
                scale only axis=true,
                width=0.5\textwidth,
                height=0.5\textwidth,
                xlabel=viscosity $\nu$,
                ylabel= $\left\Vert \uVec_{T}  -\uVec \right\Vert_{L^2}$,
                label style={anchor=near ticklabel, font=\LARGE},
                ticks=major,
                tick pos=left,
                tick align=center,
                xmin=1E-5,
                xmax=1E1,
                ymin=1E-10,
                ymax=1E-6,
                legend style={cells={anchor=west}},
                legend pos=north east,
                ]
\addplot    [
            thick,
            mark = *,
            mark size=3pt,
            mark options={black, fill=white}
            ]
            table[
            x expr = {\thisrowno{0}},
            y expr = {\thisrowno{1}},
            header=false]
            {PR_ni_vs_u_2k.dat};
\addplot    [
            thick,
            densely dashed,
            mark = triangle*,
            mark size=3pt,
            mark options={black, solid, fill=white}
            ]
            table[
            x expr = {\thisrowno{0}},
            y expr = {\thisrowno{2}},
            header=false]
            {PR_ni_vs_u_2k.dat};

\legend {\hhohdiv , \hhohll }
\end{axis}
\end{tikzpicture}
\begin{tikzpicture}[scale=0.65]
\begin{axis}    [
                unbounded coords=jump,
                ymode=log,xmode=log,
                scale only axis=true,
                width=0.5\textwidth,
                height=0.5\textwidth,
                xlabel=viscosity $\nu$,
                ylabel= $\left\Vert p_{T}  -p \right\Vert_{L^2}$,
                label style={anchor=near ticklabel, font=\LARGE},
                ticks=major,
                tick pos=left,
                tick align=center,
                xmin=1E-5,
                xmax=1E+1,
                ymin=1E-8,
                ymax=1E-5,
                legend style={cells={anchor=west}},
                legend pos=north east,
                ]
\addplot    [
            thick,
            mark = *,
            mark size=3pt,
            mark options={black, fill=white}
            ]
            table[
            x expr = {\thisrowno{0}},
            y expr = {\thisrowno{1}},
            header=false]
            {PR_ni_vs_p_2k.dat};
\addplot    [
            thick,
            densely dashed,
            mark = triangle*,
            mark size=3pt,
            mark options={black, solid, fill=white}
            ]
            table[
            x expr = {\thisrowno{0}},
            y expr = {\thisrowno{2}},
            header=false]
            {PR_ni_vs_p_2k.dat};

\legend {\hhohdiv\ , \hhohll }
\end{axis}
\end{tikzpicture}
\caption{LLMS pressure gradient. $k=3$ \hhohdiv and \hhohll formulations on a $2048$ triangular elements mesh. 
         \emph{Left and right}: Velocity and pressure errors in $L^2$-norm, respectively, while changing the viscosity ($\nu$).}\label{fig:pressure_robustness_nu_analysis}
\end{figure}
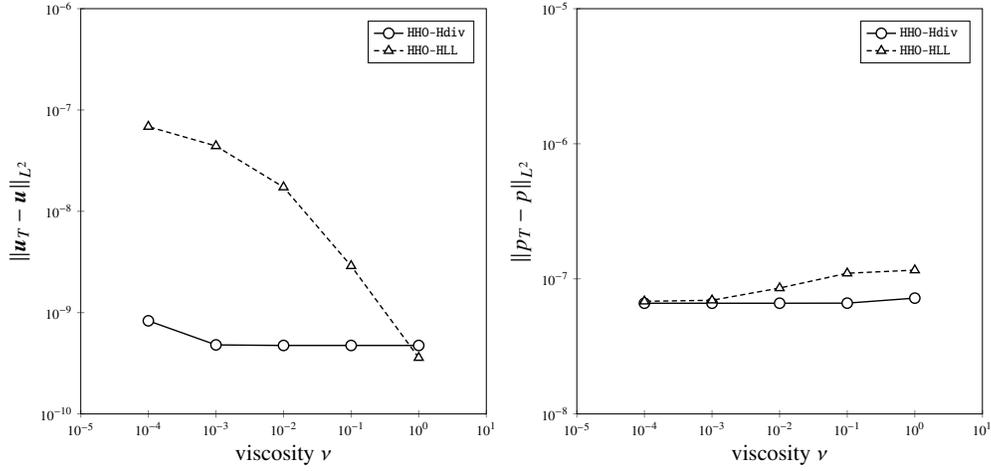

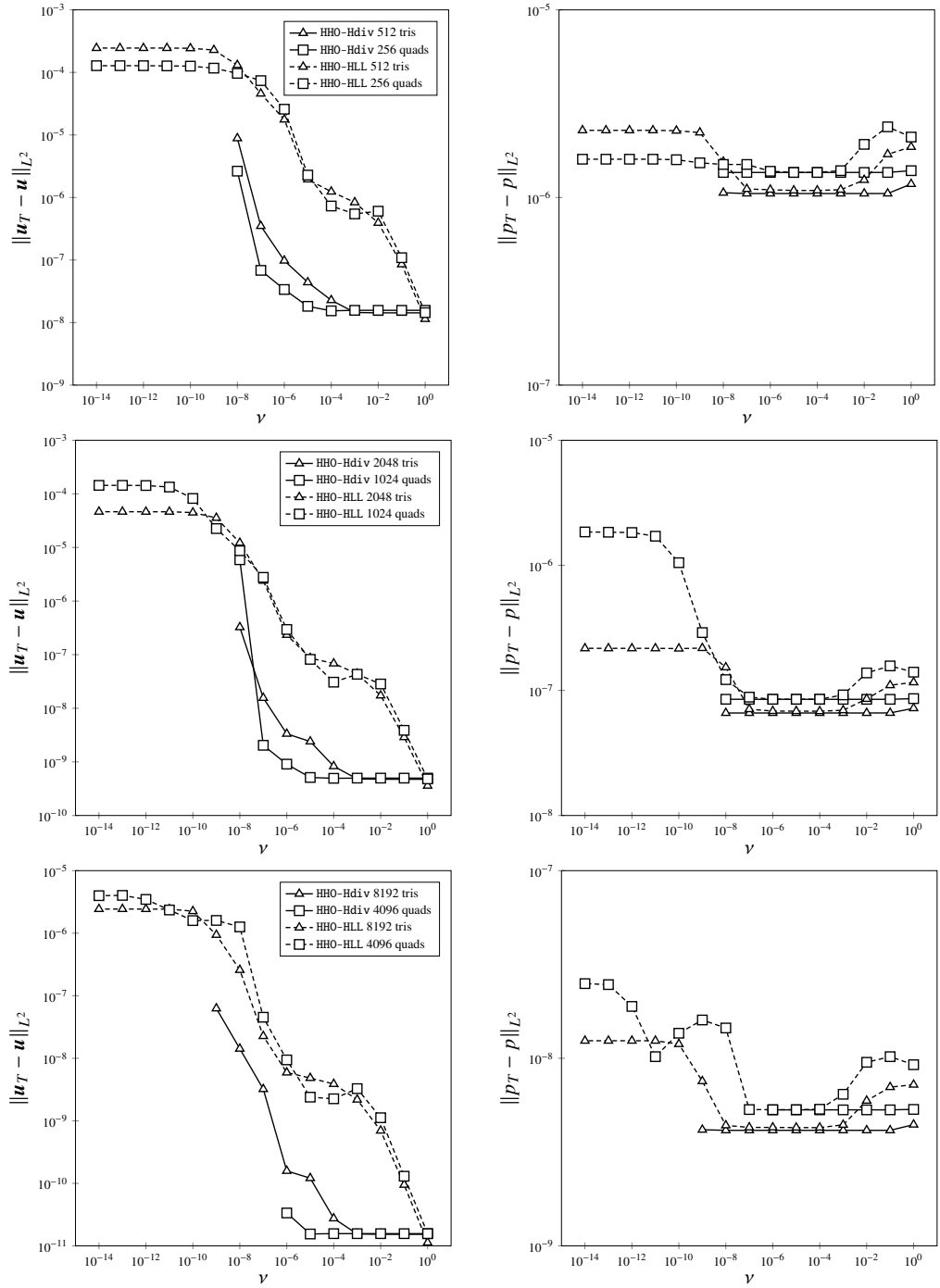
\begin{figure} [!t]
\centering
\begin{tikzpicture}[scale=0.65]
\begin{axis}    [
                unbounded coords=jump,
                ymode=log,xmode=log,
                scale only axis=true,
                width=0.5\textwidth,
                height=0.5\textwidth,
                xlabel=$\nu$,
                ylabel= $\left\Vert \uVec_{T}  -\uVec \right\Vert_{L^2}$,
                label style={anchor=near ticklabel, font=\LARGE},
                ticks=major,
                tick pos=left,
                tick align=center,
                xmin=1E-15,
                xmax=1E1,
                ymin=1E-9,
                ymax=1E-3,
                legend style={cells={anchor=west}},
                legend pos=north east,
                ]
\addplot    [
            thick,
            mark = triangle*,
            mark size=3pt,
            mark options={black, fill=white}
            ]
            table[
            x expr = {\thisrowno{0}},
            y expr = {\thisrowno{1}},
            header=false]
            {PR_ni_vs_u_256_p3.dat};
\addplot    [
            thick,
            mark = square*,
            mark size=3pt,
            mark options={black, fill=white}
            ]
            table[
            x expr = {\thisrowno{0}},
            y expr = {\thisrowno{2}},
            header=false]
            {PR_ni_vs_u_256_p3.dat};
\addplot    [
            thick,
            densely dashed,
            mark = triangle*,
            mark size=3pt,
            mark options={black, solid, fill=white}
            ]
            table[
            x expr = {\thisrowno{0}},
            y expr = {\thisrowno{3}},
            header=false]
            {PR_ni_vs_u_256_p3.dat};
\addplot    [
            thick,
            densely dashed,
            mark = square*,
            mark size=3pt,
            mark options={black, solid, fill=white}
            ]
            table[
            x expr = {\thisrowno{0}},
            y expr = {\thisrowno{4}},
            header=false]
            {PR_ni_vs_u_256_p3.dat};

\legend {\hhohdiv 512 \tris, \hhohdiv 256 \quads, \hhohll 512 \tris, \hhohll 256 \quads}
\end{axis}
\end{tikzpicture}
\hspace{5mm}
\begin{tikzpicture}[scale=0.65]
\begin{axis}    [
                unbounded coords=jump,
                ymode=log,xmode=log,
                scale only axis=true,
                width=0.5\textwidth,
                height=0.5\textwidth,
                xlabel=$\nu$,
                ylabel= $\left\Vert p_{T}  -p \right\Vert_{L^2}$,
                label style={anchor=near ticklabel, font=\LARGE},
                ticks=major,
                tick pos=left,
                tick align=center,
                xmin=1E-15,
                xmax=1E1,
                ymin=1E-7,
                ymax=1E-5,
                legend style={cells={anchor=west}},
                legend pos=north east,
                ]
\addplot    [
            thick,
            mark = triangle*,
            mark size=3pt,
            mark options={black, fill=white}
            ]
            table[
            x expr = {\thisrowno{0}},
            y expr = {\thisrowno{1}},
            header=false]
            {PR_ni_vs_p_256_p3.dat};
\addplot    [
            thick,
            mark = square*,
            mark size=3pt,
            mark options={black, fill=white}
            ]
            table[
            x expr = {\thisrowno{0}},
            y expr = {\thisrowno{2}},
            header=false]
            {PR_ni_vs_p_256_p3.dat};
\addplot    [
            thick,
            densely dashed,
            mark = triangle*,
            mark size=3pt,
            mark options={black, solid, fill=white}
            ]
            table[
            x expr = {\thisrowno{0}},
            y expr = {\thisrowno{3}},
            header=false]
            {PR_ni_vs_p_256_p3.dat};
\addplot    [
            thick,
            densely dashed,
            mark = square*,
            mark size=3pt,
            mark options={black, solid, fill=white}
            ]
            table[
            x expr = {\thisrowno{0}},
            y expr = {\thisrowno{4}},
            header=false]
            {PR_ni_vs_p_256_p3.dat};

\end{axis}
\end{tikzpicture}
\begin{tikzpicture}[scale=0.65]
\begin{axis}    [
                unbounded coords=jump,
                ymode=log,xmode=log,
                scale only axis=true,
                width=0.5\textwidth,
                height=0.5\textwidth,
                xlabel=$\nu$,
                ylabel= $\left\Vert \uVec_{T}  -\uVec \right\Vert_{L^2}$,
                label style={anchor=near ticklabel, font=\LARGE},
                ticks=major,
                tick pos=left,
                tick align=center,
                xmin=1E-15,
                xmax=1E1,
                ymin=1E-10,
                ymax=1E-3,
                legend style={cells={anchor=west}},
                legend pos=north east,
                ]
\addplot    [
            thick,
            mark = triangle*,
            mark size=3pt,
            mark options={black, fill=white}
            ]
            table[
            x expr = {\thisrowno{0}},
            y expr = {\thisrowno{1}},
            header=false]
            {PR_ni_vs_u_1k_p3.dat};
\addplot    [
            thick,
            mark = square*,
            mark size=3pt,
            mark options={black, fill=white}
            ]
            table[
            x expr = {\thisrowno{0}},
            y expr = {\thisrowno{2}},
            header=false]
            {PR_ni_vs_u_1k_p3.dat};
\addplot    [
            thick,
            densely dashed,
            mark = triangle*,
            mark size=3pt,
            mark options={black, solid, fill=white}
            ]
            table[
            x expr = {\thisrowno{0}},
            y expr = {\thisrowno{3}},
            header=false]
            {PR_ni_vs_u_1k_p3.dat};
\addplot    [
            thick,
            densely dashed,
            mark = square*,
            mark size=3pt,
            mark options={black, solid, fill=white}
            ]
            table[
            x expr = {\thisrowno{0}},
            y expr = {\thisrowno{4}},
            header=false]
            {PR_ni_vs_u_1k_p3.dat};

\legend {\hhohdiv 2048 \tris, \hhohdiv 1024 \quads, \hhohll 2048 \tris, \hhohll 1024 \quads}
\end{axis}
\end{tikzpicture}
\hspace{5mm}
\begin{tikzpicture}[scale=0.65]
\begin{axis}    [
                unbounded coords=jump,
                ymode=log,xmode=log,
                scale only axis=true,
                width=0.5\textwidth,
                height=0.5\textwidth,
                xlabel=$\nu$,
                ylabel= $\left\Vert p_{T}  -p \right\Vert_{L^2}$,
                label style={anchor=near ticklabel, font=\LARGE},
                ticks=major,
                tick pos=left,
                tick align=center,
                xmin=1E-15,
                xmax=1E1,
                ymin=1E-8,
                ymax=1E-5,
                legend style={cells={anchor=west}},
                legend pos=north east,
                ]
\addplot    [
            thick,
            mark = triangle*,
            mark size=3pt,
            mark options={black, fill=white}
            ]
            table[
            x expr = {\thisrowno{0}},
            y expr = {\thisrowno{1}},
            header=false]
            {PR_ni_vs_p_1k_p3.dat};
\addplot    [
            thick,
            mark = square*,
            mark size=3pt,
            mark options={black, fill=white}
            ]
            table[
            x expr = {\thisrowno{0}},
            y expr = {\thisrowno{2}},
            header=false]
            {PR_ni_vs_p_1k_p3.dat};
\addplot    [
            thick,
            densely dashed,
            mark = triangle*,
            mark size=3pt,
            mark options={black, solid, fill=white}
            ]
            table[
            x expr = {\thisrowno{0}},
            y expr = {\thisrowno{3}},
            header=false]
            {PR_ni_vs_p_1k_p3.dat};
\addplot    [
            thick,
            densely dashed,
            mark = square*,
            mark size=3pt,
            mark options={black, solid, fill=white}
            ]
            table[
            x expr = {\thisrowno{0}},
            y expr = {\thisrowno{4}},
            header=false]
            {PR_ni_vs_p_1k_p3.dat};

\end{axis}
\end{tikzpicture}
\begin{tikzpicture}[scale=0.65]
\begin{axis}    [
                unbounded coords=jump,
                ymode=log,xmode=log,
                scale only axis=true,
                width=0.5\textwidth,
                height=0.5\textwidth,
                xlabel=$\nu$,
                ylabel= $\left\Vert \uVec_{T}  -\uVec \right\Vert_{L^2}$,
                label style={anchor=near ticklabel, font=\LARGE},
                ticks=major,
                tick pos=left,
                tick align=center,
                xmin=1E-15,
                xmax=1E1,
                ymin=1E-11,
                ymax=1E-5,
                legend style={cells={anchor=west}},
                legend pos=north east,
                ]
\addplot    [
            thick,
            mark = triangle*,
            mark size=3pt,
            mark options={black, fill=white}
            ]
            table[
            x expr = {\thisrowno{0}},
            y expr = {\thisrowno{1}},
            header=false]
            {PR_ni_vs_u_4k_p3.dat};
\addplot    [
            thick,
            mark = square*,
            mark size=3pt,
            mark options={black, fill=white}
            ]
            table[
            x expr = {\thisrowno{0}},
            y expr = {\thisrowno{2}},
            header=false]
            {PR_ni_vs_u_4k_p3.dat};
\addplot    [
            thick,
            densely dashed,
            mark = triangle*,
            mark size=3pt,
            mark options={black, solid, fill=white}
            ]
            table[
            x expr = {\thisrowno{0}},
            y expr = {\thisrowno{3}},
            header=false]
            {PR_ni_vs_u_4k_p3.dat};
\addplot    [
            thick,
            densely dashed,
            mark = square*,
            mark size=3pt,
            mark options={black, solid, fill=white}
            ]
            table[
            x expr = {\thisrowno{0}},
            y expr = {\thisrowno{4}},
            header=false]
            {PR_ni_vs_u_4k_p3.dat};

\legend {\hhohdiv 8192 \tris, \hhohdiv 4096 \quads, \hhohll 8192 \tris, \hhohll 4096 \quads}
\end{axis}
\end{tikzpicture}
\hspace{5mm}
\begin{tikzpicture}[scale=0.65]
\begin{axis}    [
                unbounded coords=jump,
                ymode=log,xmode=log,
                scale only axis=true,
                width=0.5\textwidth,
                height=0.5\textwidth,
                xlabel=$\nu$,
                ylabel= $\left\Vert p_{T}  -p \right\Vert_{L^2}$,
                label style={anchor=near ticklabel, font=\LARGE},
                ticks=major,
                tick pos=left,
                tick align=center,
                xmin=1E-15,
                xmax=1E1,
                ymin=1E-9,
                ymax=1E-7,
                legend style={cells={anchor=west}},
                legend pos=north east,
                ]
\addplot    [
            thick,
            mark = triangle*,
            mark size=3pt,
            mark options={black, fill=white}
            ]
            table[
            x expr = {\thisrowno{0}},
            y expr = {\thisrowno{1}},
            header=false]
            {PR_ni_vs_p_4k_p3.dat};
\addplot    [
            thick,
            mark = square*,
            mark size=3pt,
            mark options={black, fill=white}
            ]
            table[
            x expr = {\thisrowno{0}},
            y expr = {\thisrowno{2}},
            header=false]
            {PR_ni_vs_p_4k_p3.dat};
\addplot    [
            thick,
            densely dashed,
            mark = triangle*,
            mark size=3pt,
            mark options={black, solid, fill=white}
            ]
            table[
            x expr = {\thisrowno{0}},
            y expr = {\thisrowno{3}},
            header=false]
            {PR_ni_vs_p_4k_p3.dat};
\addplot    [
            thick,
            densely dashed,
            mark = square*,
            mark size=3pt,
            mark options={black, solid, fill=white}
            ]
            table[
            x expr = {\thisrowno{0}},
            y expr = {\thisrowno{4}},
            header=false]
            {PR_ni_vs_p_4k_p3.dat};

\end{axis}
\end{tikzpicture}
\caption{LLMS pressure gradient. Robustness in the inviscid limit, $k=3$ \hhohdiv (solid lines) and \hhohll (dashed lines) formulations. \emph{Left and right}: velocity and pressure error in $L^2$-norm, respectively. 
         In each row, from top to bottom, increasingly dense triangular (triangular marks) and quadrilateral (square marks) meshes are considered.}\label{fig:pressure_robustness_inviscid_p3}
\end{figure}
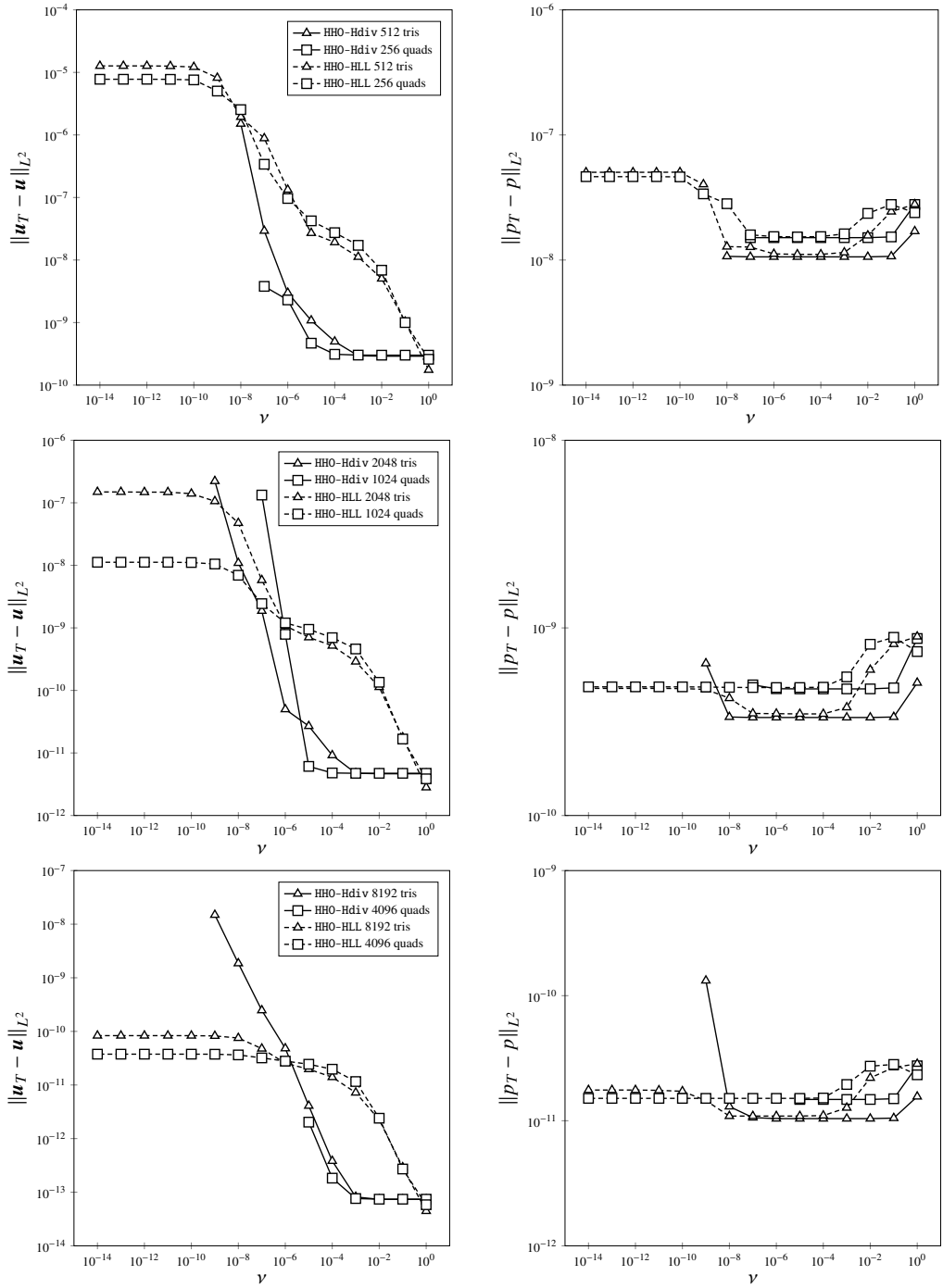
\begin{figure} [!t]
\centering
\begin{tikzpicture}[scale=0.65]
\begin{axis}    [
                unbounded coords=jump,
                ymode=log,xmode=log,
                scale only axis=true,
                width=0.5\textwidth,
                height=0.5\textwidth,
                xlabel=$\nu$,
                ylabel= $\left\Vert \uVec_{T}  -\uVec \right\Vert_{L^2}$,
                label style={anchor=near ticklabel, font=\LARGE},
                ticks=major,
                tick pos=left,
                tick align=center,
                xmin=1E-15,
                xmax=1E1,
                ymin=1E-10,
                ymax=1E-4,
                legend style={cells={anchor=west}},
                legend pos=north east,
                ]
\addplot    [
            thick,
            mark = triangle*,
            mark size=3pt,
            mark options={black, fill=white}
            ]
            table[
            x expr = {\thisrowno{0}},
            y expr = {\thisrowno{1}},
            header=false]
            {PR_ni_vs_u_256_p4.dat};
\addplot    [
            thick,
            mark = square*,
            mark size=3pt,
            mark options={black, fill=white}
            ]
            table[
            x expr = {\thisrowno{0}},
            y expr = {\thisrowno{2}},
            header=false]
            {PR_ni_vs_u_256_p4.dat};
\addplot    [
            thick,
            densely dashed,
            mark = triangle*,
            mark size=3pt,
            mark options={black, solid, fill=white}
            ]
            table[
            x expr = {\thisrowno{0}},
            y expr = {\thisrowno{3}},
            header=false]
            {PR_ni_vs_u_256_p4.dat};
\addplot    [
            thick,
            densely dashed,
            mark = square*,
            mark size=3pt,
            mark options={black, solid, fill=white}
            ]
            table[
            x expr = {\thisrowno{0}},
            y expr = {\thisrowno{4}},
            header=false]
            {PR_ni_vs_u_256_p4.dat};

\legend {\hhohdiv 512 \tris, \hhohdiv 256 \quads, \hhohll 512 \tris, \hhohll 256 \quads}
\end{axis}
\end{tikzpicture}
\hspace{5mm}
\begin{tikzpicture}[scale=0.65]
\begin{axis}    [
                unbounded coords=jump,
                ymode=log,xmode=log,
                scale only axis=true,
                width=0.5\textwidth,
                height=0.5\textwidth,
                xlabel=$\nu$,
                ylabel= $\left\Vert p_{T}  -p \right\Vert_{L^2}$,
                label style={anchor=near ticklabel, font=\LARGE},
                ticks=major,
                tick pos=left,
                tick align=center,
                xmin=1E-15,
                xmax=1E1,
                ymin=1E-9,
                ymax=1E-6,
                legend style={cells={anchor=west}},
                legend pos=north east,
                ]
\addplot    [
            thick,
            mark = triangle*,
            mark size=3pt,
            mark options={black, fill=white}
            ]
            table[
            x expr = {\thisrowno{0}},
            y expr = {\thisrowno{1}},
            header=false]
            {PR_ni_vs_p_256_p4.dat};
\addplot    [
            thick,
            mark = square*,
            mark size=3pt,
            mark options={black, fill=white}
            ]
            table[
            x expr = {\thisrowno{0}},
            y expr = {\thisrowno{2}},
            header=false]
            {PR_ni_vs_p_256_p4.dat};
\addplot    [
            thick,
            densely dashed,
            mark = triangle*,
            mark size=3pt,
            mark options={black, solid, fill=white}
            ]
            table[
            x expr = {\thisrowno{0}},
            y expr = {\thisrowno{3}},
            header=false]
            {PR_ni_vs_p_256_p4.dat};
\addplot    [
            thick,
            densely dashed,
            mark = square*,
            mark size=3pt,
            mark options={black, solid, fill=white}
            ]
            table[
            x expr = {\thisrowno{0}},
            y expr = {\thisrowno{4}},
            header=false]
            {PR_ni_vs_p_256_p4.dat};

\end{axis}
\end{tikzpicture}
\begin{tikzpicture}[scale=0.65]
\begin{axis}    [
                unbounded coords=jump,
                ymode=log,xmode=log,
                scale only axis=true,
                width=0.5\textwidth,
                height=0.5\textwidth,
                xlabel=$\nu$,
                ylabel= $\left\Vert \uVec_{T}  -\uVec \right\Vert_{L^2}$,
                label style={anchor=near ticklabel, font=\LARGE},
                ticks=major,
                tick pos=left,
                tick align=center,
                xmin=1E-15,
                xmax=1E1,
                ymin=1E-12,
                ymax=1E-6,
                legend style={cells={anchor=west}},
                legend pos=north east,
                ]
\addplot    [
            thick,
            mark = triangle*,
            mark size=3pt,
            mark options={black, fill=white}
            ]
            table[
            x expr = {\thisrowno{0}},
            y expr = {\thisrowno{1}},
            header=false]
            {PR_ni_vs_u_1k_p4.dat};
\addplot    [
            thick,
            mark = square*,
            mark size=3pt,
            mark options={black, fill=white}
            ]
            table[
            x expr = {\thisrowno{0}},
            y expr = {\thisrowno{2}},
            header=false]
            {PR_ni_vs_u_1k_p4.dat};
\addplot    [
            thick,
            densely dashed,
            mark = triangle*,
            mark size=3pt,
            mark options={black, solid, fill=white}
            ]
            table[
            x expr = {\thisrowno{0}},
            y expr = {\thisrowno{3}},
            header=false]
            {PR_ni_vs_u_1k_p4.dat};
\addplot    [
            thick,
            densely dashed,
            mark = square*,
            mark size=3pt,
            mark options={black, solid, fill=white}
            ]
            table[
            x expr = {\thisrowno{0}},
            y expr = {\thisrowno{4}},
            header=false]
            {PR_ni_vs_u_1k_p4.dat};

\legend {\hhohdiv 2048 \tris, \hhohdiv 1024 \quads, \hhohll 2048 \tris, \hhohll 1024 \quads}
\end{axis}
\end{tikzpicture}
\hspace{5mm}
\begin{tikzpicture}[scale=0.65]
\begin{axis}    [
                unbounded coords=jump,
                ymode=log,xmode=log,
                scale only axis=true,
                width=0.5\textwidth,
                height=0.5\textwidth,
                xlabel=$\nu$,
                ylabel= $\left\Vert p_{T}  -p \right\Vert_{L^2}$,
                label style={anchor=near ticklabel, font=\LARGE},
                ticks=major,
                tick pos=left,
                tick align=center,
                xmin=1E-15,
                xmax=1E1,
                ymin=1E-10,
                ymax=1E-8,
                legend style={cells={anchor=west}},
                legend pos=north east,
                ]
\addplot    [
            thick,
            mark = triangle*,
            mark size=3pt,
            mark options={black, fill=white}
            ]
            table[
            x expr = {\thisrowno{0}},
            y expr = {\thisrowno{1}},
            header=false]
            {PR_ni_vs_p_1k_p4.dat};
\addplot    [
            thick,
            mark = square*,
            mark size=3pt,
            mark options={black, fill=white}
            ]
            table[
            x expr = {\thisrowno{0}},
            y expr = {\thisrowno{2}},
            header=false]
            {PR_ni_vs_p_1k_p4.dat};
\addplot    [
            thick,
            densely dashed,
            mark = triangle*,
            mark size=3pt,
            mark options={black, solid, fill=white}
            ]
            table[
            x expr = {\thisrowno{0}},
            y expr = {\thisrowno{3}},
            header=false]
            {PR_ni_vs_p_1k_p4.dat};
\addplot    [
            thick,
            densely dashed,
            mark = square*,
            mark size=3pt,
            mark options={black, solid, fill=white}
            ]
            table[
            x expr = {\thisrowno{0}},
            y expr = {\thisrowno{4}},
            header=false]
            {PR_ni_vs_p_1k_p4.dat};

\end{axis}
\end{tikzpicture}
\begin{tikzpicture}[scale=0.65]
\begin{axis}    [
                unbounded coords=jump,
                ymode=log,xmode=log,
                scale only axis=true,
                width=0.5\textwidth,
                height=0.5\textwidth,
                xlabel=$\nu$,
                ylabel= $\left\Vert \uVec_{T}  -\uVec \right\Vert_{L^2}$,
                label style={anchor=near ticklabel, font=\LARGE},
                ticks=major,
                tick pos=left,
                tick align=center,
                xmin=1E-15,
                xmax=1E1,
                ymin=1E-14,
                ymax=1E-7,
                legend style={cells={anchor=west}},
                legend pos=north east,
                ]
\addplot    [
            thick,
            mark = triangle*,
            mark size=3pt,
            mark options={black, fill=white}
            ]
            table[
            x expr = {\thisrowno{0}},
            y expr = {\thisrowno{1}},
            header=false]
            {PR_ni_vs_u_4k_p4.dat};
\addplot    [
            thick,
            mark = square*,
            mark size=3pt,
            mark options={black, fill=white}
            ]
            table[
            x expr = {\thisrowno{0}},
            y expr = {\thisrowno{2}},
            header=false]
            {PR_ni_vs_u_4k_p4.dat};
\addplot    [
            thick,
            densely dashed,
            mark = triangle*,
            mark size=3pt,
            mark options={black, solid, fill=white}
            ]
            table[
            x expr = {\thisrowno{0}},
            y expr = {\thisrowno{3}},
            header=false]
            {PR_ni_vs_u_4k_p4.dat};
\addplot    [
            thick,
            densely dashed,
            mark = square*,
            mark size=3pt,
            mark options={black, solid, fill=white}
            ]
            table[
            x expr = {\thisrowno{0}},
            y expr = {\thisrowno{4}},
            header=false]
            {PR_ni_vs_u_4k_p4.dat};

\legend {\hhohdiv 8192 \tris, \hhohdiv 4096 \quads, \hhohll 8192 \tris, \hhohll 4096 \quads}
\end{axis}
\end{tikzpicture}
\hspace{5mm}
\begin{tikzpicture}[scale=0.65]
\begin{axis}    [
                unbounded coords=jump,
                ymode=log,xmode=log,
                scale only axis=true,
                width=0.5\textwidth,
                height=0.5\textwidth,
                xlabel=$\nu$,
                ylabel= $\left\Vert p_{T}  -p \right\Vert_{L^2}$,
                label style={anchor=near ticklabel, font=\LARGE},
                ticks=major,
                tick pos=left,
                tick align=center,
                xmin=1E-15,
                xmax=1E1,
                ymin=1E-12,
                ymax=1E-9,
                legend style={cells={anchor=west}},
                legend pos=north east,
                ]
\addplot    [
            thick,
            mark = triangle*,
            mark size=3pt,
            mark options={black, fill=white}
            ]
            table[
            x expr = {\thisrowno{0}},
            y expr = {\thisrowno{1}},
            header=false]
            {PR_ni_vs_p_4k_p4.dat};
\addplot    [
            thick,
            mark = square*,
            mark size=3pt,
            mark options={black, fill=white}
            ]
            table[
            x expr = {\thisrowno{0}},
            y expr = {\thisrowno{2}},
            header=false]
            {PR_ni_vs_p_4k_p4.dat};
\addplot    [
            thick,
            densely dashed,
            mark = triangle*,
            mark size=3pt,
            mark options={black, solid, fill=white}
            ]
            table[
            x expr = {\thisrowno{0}},
            y expr = {\thisrowno{3}},
            header=false]
            {PR_ni_vs_p_4k_p4.dat};
\addplot    [
            thick,
            densely dashed,
            mark = square*,
            mark size=3pt,
            mark options={black, solid, fill=white}
            ]
            table[
            x expr = {\thisrowno{0}},
            y expr = {\thisrowno{4}},
            header=false]
            {PR_ni_vs_p_4k_p4.dat};

\end{axis}
\end{tikzpicture}
\caption{LLMS pressure gradient. Robustness in the inviscid limit, $k=4$ \hhohdiv (solid lines) and \hhohll (dashed lines) formulations. \emph{Left and right}: velocity and pressure error in $L^2$-norm, respectively. 
         In each row, from top to bottom, increasingly dense triangular (triangular marks) and quadrilateral (square marks) meshes are considered.}\label{fig:pressure_robustness_inviscid_p4}
\end{figure}

\subsection{Gresho-Chan vortex}
\label{sec:GCVortex}
Originally proposed by Gresho \& Chan \cite{Gresho.Chan:90}, who named it `triangular vortex', 
this viscous 2D model problem is designed to study the schemes capability of preserving vortical structures. 
Our analysis relies on the numerical set-up described by Gauger \ea \cite{Gauger.Linke.Schroeder:19}: time integration is carried out 
over a double periodic unit squared domain $\Omega=\left(0,\,1\right)\times\left(0,\,1\right)$
based on the following initial conditions 
\begin{equation*}
	\left(u,\,v,\,p\right)
		=	\left(
			u_0,\,v_0,\,p_0
			\right)
			+\begin{cases}
            	\bigg(
                -5\widetilde{y},\,   5\widetilde{x},\,	\dfrac{25}{2}r^{2}
                \bigg),
                &  \mbox{if} \; 0\leq r<0.2,   \\
            	\bigg(
               	-2\dfrac{\widetilde{y}}{r}  +5\widetilde{y},\, 2\dfrac{\widetilde{x}}{r}  -5\widetilde{x},\, 4\ln{5r}+\dfrac{25}{2}r^{2}  -20r  +4
                \bigg),
               	&  \mbox{if} \; 0.2\leq r<0.4   \\
            	\bigg(
               	0,\,	0,\,	4\ln{2}  -2
                \bigg),
                &  \mbox{if} \; 0.4\leq r;
         	\end{cases},
\end{equation*}
where $u_{0}$, $v_{0}$, $p_{0}$ are user-defined real parameters.
Notice that a local coordinate system $\widetilde{\xVec}=\xVec-\boldsymbol{x}_0$ is employed to define the radial distance 
$r=(\widetilde{x}^{2}+\widetilde{y}^{2})^{0.5}$ from the vortex center $\boldsymbol{x}_0=\left[0.5\,\,0.5\right]^{\intercal}$ at the initial simulation time.

Figure~\ref{fig:gresho-chan_radial_initial} shows the initial circumferential velocity $u_{\theta}$ 
and vorticity $\boldsymbol{\omega}=\nabla\times\uVec$ fields as function of the radial distance from the vortex center in case of a `standing vortex problem', i.e., $u_{0} =v_{0}=0$.
Notice the peculiar triangular distribution of the circumferential velocity, from whom the test name originates.
Interestingly, the velocity field yields a discontinuous vorticity field: 
a rigid body rotation core ($r<0.2$) with constant counter-clockwise vorticity comes in contact with an annular region 
($0.2\leq r<0.4$) with radially increasing clockwise vorticity followed by an external region ($0.4\leq r$) with fluid at rest (null vorticity).
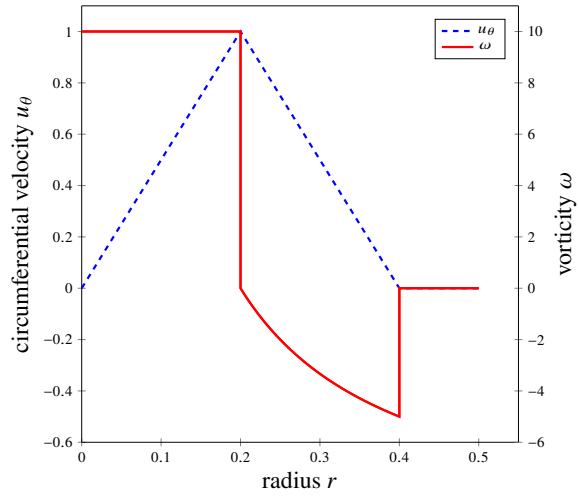
\begin{figure} [!t]
\centering
\begin{tikzpicture}[scale=0.7]
\begin{axis}    [
                unbounded coords=jump,
                scale only axis=true,
                width=0.5\textwidth,
                height=0.5\textwidth,
                axis y line*=left,
                xlabel=radius $r$,
                ylabel=circumferential velocity $u_{\theta}$,
                label style={anchor=near ticklabel, font=\LARGE},
                ticks=major,
                tick pos=left,
                tick align=center,
                xmin=0,
                xmax=0.55,
                ymin=-0.6,
                ymax=1.1,
                legend style={cells={anchor=west}},
                legend pos=north east,
                ]
\addplot    [
            blue,
            very thick,
            dashed
            ]
            table[
            x expr = {\thisrowno{0}},
            y expr = {\thisrowno{1}},
            header=false]
            {GC_initial.dat};
\addplot    [
            red,
            very thick
            ]
            table[
            x expr = {\thisrowno{0}},
            y expr = {0.1*\thisrowno{2}},
            header=false]
            {GC_initial.dat};
\legend {$u_{\theta}$, $\omega$}
\end{axis}
\begin{axis}    [
                unbounded coords=jump,
                scale only axis=true,
                width=0.5\textwidth,
                height=0.5\textwidth,
                axis y line*=right,
                axis x line=none,
                ylabel=vorticity $\omega$,
                label style={anchor=near ticklabel, font=\LARGE},
                ticks=major,
                tick align=center,
                xmin=0,
                xmax=0.55,
                ymin=-6,
                ymax=11,
                legend style={cells={anchor=west}},
                legend pos=north east,
                ]
\addplot    [
            red,
            very thick
            ]
            table[
            x expr = {\thisrowno{0}},
            y expr = {\thisrowno{2}},
            header=false]
            {GC_initial.dat};
\end{axis}
\end{tikzpicture}
    \caption{Gresho-Chan vortex. Radial distribution of the initial circumferential velocity and vorticity fields for the `standing vortex problem' ($u_{0}=v_{0}=0$).}
	\label{fig:gresho-chan_radial_initial}
\end{figure}

The numerical investigation is performed imposing $u_{0}= v_{0} = {1}/{3}$ (moving vortex), $p_{0}=0$ and setting the viscosity as $\nu=10^{-5}$.
BDF2 implicit time integration \cite{Curtiss.Hirschfelder:1952} with a constant step size $\Delta t = t_{F}/10^{4}$ is adopted to advance the solution 
up to the dimensionless end time $t_{F}=3$, resulting in the simulation of one period of the moving vortex.
At each time step the nonlinear system is solved using a Newton-Krylov iterative method \cite{Brown.Saad:94}.
We consider $k=4$ and $k=7$ \hhohdiv and \hhohll formulations on a regular $32 \times 32$ quadrilateral elements mesh.
The higher polynomial degree is considered as a reference to evaluate the behavior of the lowest order discretization 
with respect to the evolution of kinetic energy $\mathcal{K}$ and enstrophy $\mathcal{E}$, defined as follows
\begin{equation*}
	\mathcal{K}	=	\dfrac{1}{2}\int_{\Omega} \uVec \cdot \uVec,
	\qquad
	\mathcal{E}	=	\dfrac{1}{2}\int_{\Omega} \boldsymbol{\omega}\cdot\boldsymbol{\omega}.
\end{equation*}

All the simulations performed provide an accurate picture of kinetic energy evolution, see Figure~\ref{fig:gresho-chan_kinetic_enstrophy}.
Notice in particular that the kinetic energy decays are superimposed for both schemes and both polynomial degrees. 
Some differences might be appreciated in terms of enstrophy behavior: 
while the \hhohll and \hhohdiv curves are almost superimposed at $k=7$, at $k=4$ the \hhohll scheme underestimates the enstrophy by a more significant amount.   %
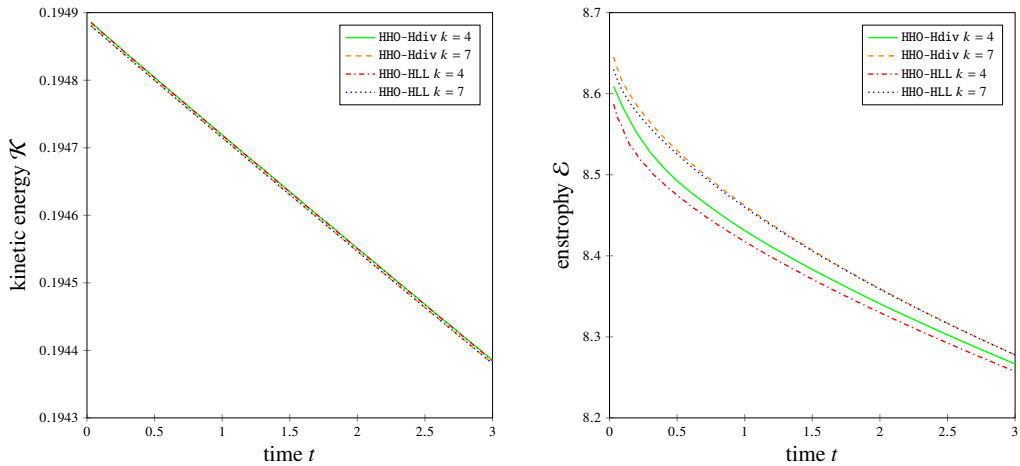
\begin{figure} [!t]
\begin{center}
\begin{tikzpicture}[scale=0.65]
\begin{axis}    [
                unbounded coords=jump,
                scale only axis=true,
                width=0.5\textwidth,
                height=0.5\textwidth,
                xlabel=time $t$,
                ylabel=kinetic energy $\mathcal{K}$,
                label style={anchor=near ticklabel, font=\LARGE},
                y tick label style={/pgf/number format/.cd, fixed, precision=6},
                ticks=major,
                tick pos=left,
                tick align=center,
                xmin=0,
                xmax=3,
                ymin=0.1943,
                ymax=0.1949,
                legend style={cells={anchor=west, font=\footnotesize}},
                legend pos=north east,
                ]
\addplot    [
            green,
            thick,
            ]
            table[
            x expr = {\thisrowno{0}},
            y expr = {\thisrowno{1}},
            header=false]
            {GC-Hdiv_p4.dat};
\addplot    [
            orange,
            thick,
            dashed,
            ]
            table[
            x expr = {\thisrowno{0}},
            y expr = {\thisrowno{1}},
            header=false]
            {GC-Hdiv_p7.dat};
\addplot    [
            red,
            thick,
            dashdotted,
            ]
            table[
            x expr = {\thisrowno{0}},
            y expr = {\thisrowno{1}},
            header=false]
            {GC-HLL_p4.dat}; 
\addplot    [
            blue,
            thick,
            dotted,
            ]
            table[
            x expr = {\thisrowno{0}},
            y expr = {\thisrowno{1}},
            header=false]
            {GC-HLL_p7.dat};

\legend {\hhohdiv $k=4$,\hhohdiv $k=7$, \hhohll $k=4$, \hhohll $k=7$}
\end{axis}
\end{tikzpicture}
\hspace{5mm}
\begin{tikzpicture}[scale=0.65]
\begin{axis}    [
                unbounded coords=jump,
                scale only axis=true,
                width=0.5\textwidth,
                height=0.5\textwidth,
                xlabel=time $t$,
                ylabel=enstrophy $\mathcal{E}$,
                label style={anchor=near ticklabel, font=\LARGE},
                y tick label style={/pgf/number format/.cd, fixed, precision=6},
                ticks=major,
                tick pos=left,
                tick align=center,
                xmin=0,
                xmax=3,
                ymin=8.2,
                ymax=8.7,
                legend style={cells={anchor=west, font=\footnotesize}},
                legend pos=north east,
                ]
\addplot    [
            green,
            thick,
            ]
            table[
            x expr = {\thisrowno{0}},
            y expr = {\thisrowno{2}},
            header=false]
            {GC-Hdiv_p4.dat};
\addplot    [
            orange,
            thick,
            dashed,
            ]
            table[
            x expr = {\thisrowno{0}},
            y expr = {\thisrowno{2}},
            header=false]
            {GC-Hdiv_p7.dat};
\addplot    [
            red,
            thick,
            dashdotted,
            ]
            table[
            x expr = {\thisrowno{0}},
            y expr = {\thisrowno{2}},
            header=false]
            {GC-HLL_p4.dat}; 
\addplot    [
            blue,
            thick,
            dotted,
            ]
            table[
            x expr = {\thisrowno{0}},
            y expr = {\thisrowno{2}},
            header=false]
            {GC-HLL_p7.dat};

\legend {\hhohdiv $k=4$,\hhohdiv $k=7$, \hhohll $k=4$, \hhohll $k=7$}
\end{axis}
\end{tikzpicture}
\end{center}
    \caption{Gresho-Chan vortex. Kinetic energy $\mathcal{K}$ and enstrophy $\mathcal{E}$ time evolution.}
\label{fig:gresho-chan_kinetic_enstrophy}
\end{figure}

Figure~\ref{fig:gresho-chan_vorticity} depicts the vorticity fields at the initial and final simulation times.
Because of the discontinuous nature of the initial vorticity field, local over/under-shoots are present at initialization, 
nonetheless, the final solution looks smoother and the vortex is well resolved.
The results suggest that both schemes are capable of satisfactorily preserving vortical structures.
Note however that \hhohdiv shows some high-frequency oscillations that are absent in the final \hhohdiv solution.
\begin{figure}[!t]
    \centering
	\includegraphics[width=0.4\textwidth]{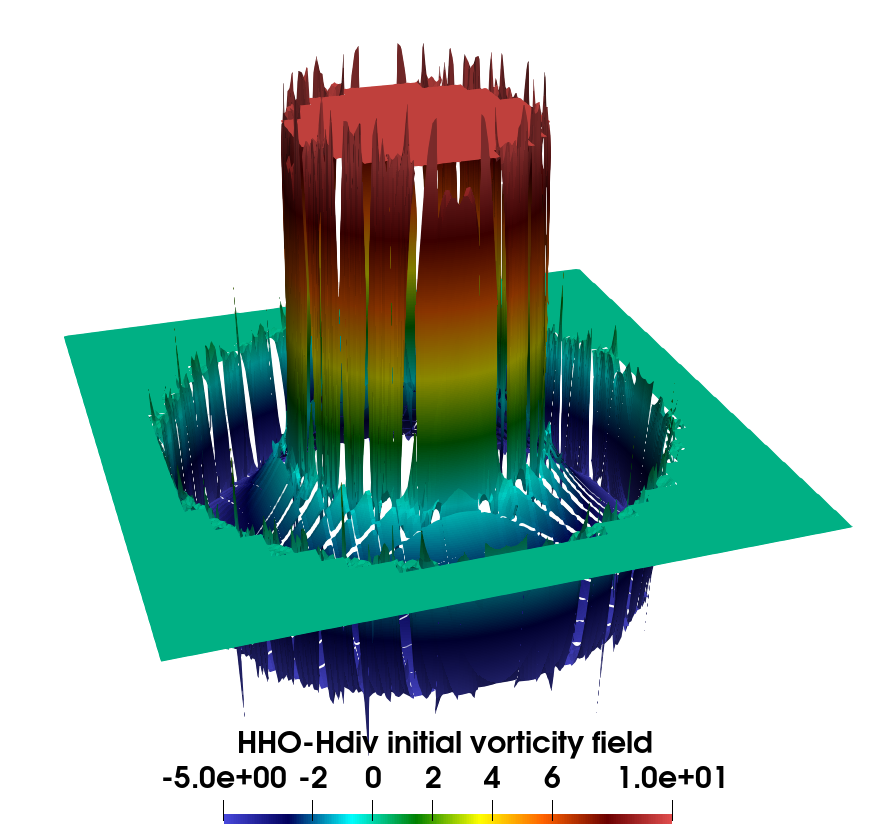}
	\includegraphics[width=0.4\textwidth]{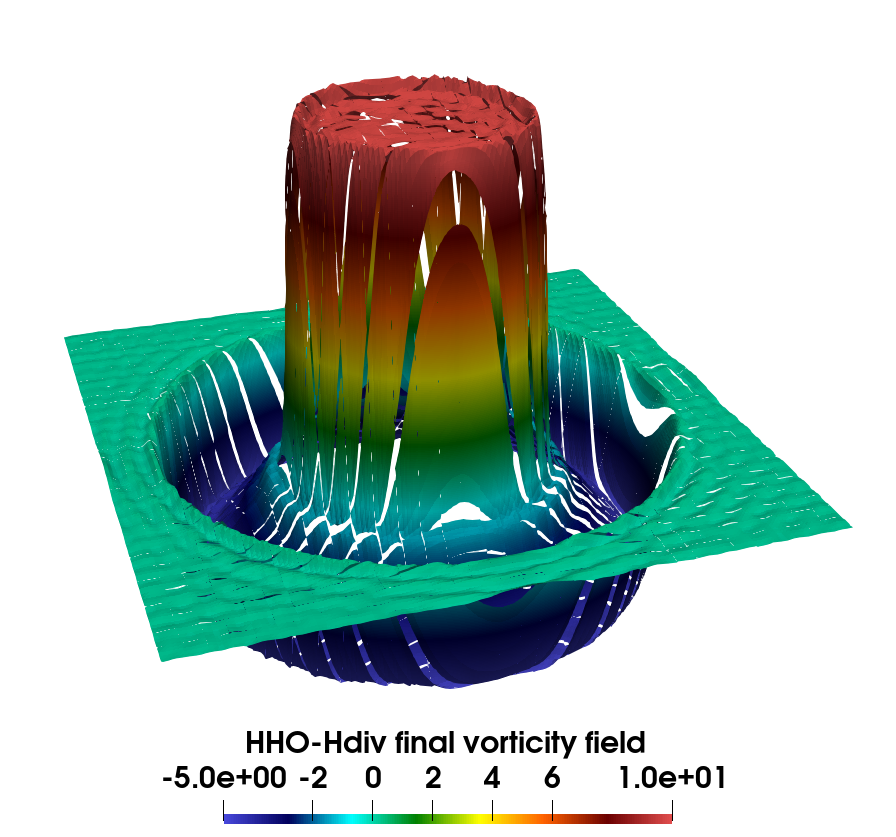}
	\includegraphics[width=0.4\textwidth]{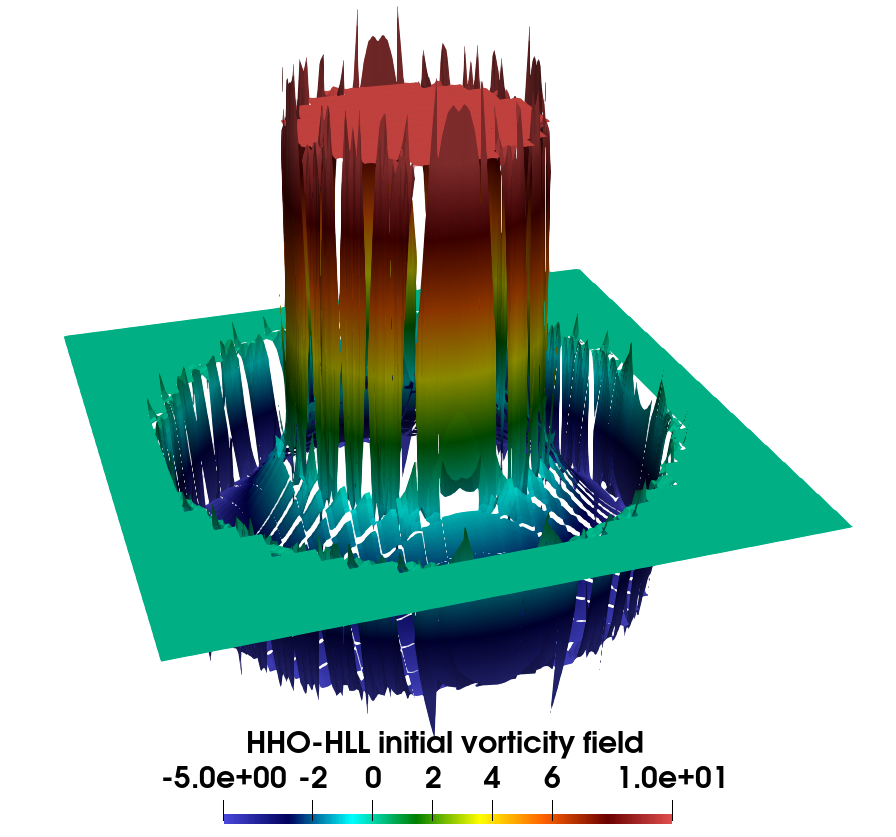}
	\includegraphics[width=0.4\textwidth]{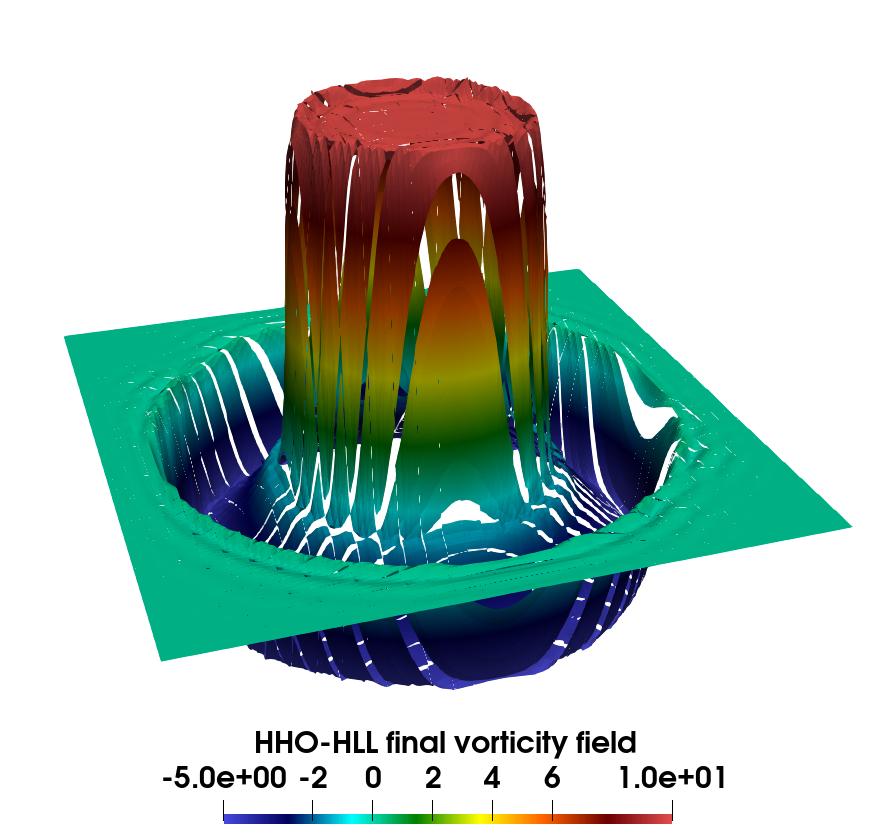}
    \caption{Gresho-Chan vortex. Vorticity fields (color coded and warped by vorticity magnitude) at $k=4$. 
             \emph{Left and right}: initial and final simulation times, respectively. \emph{Top and bottom}: \hhohdiv and \hhohll schemes, respectively.}
	\label{fig:gresho-chan_vorticity}
\end{figure}

\subsection{Double shear layer}
\label{sec:DSLayer}
The double shear layer problem devised by Bell \ea \cite{Bell.et.al:1989} focuses on inviscid ($\nu=0$) unsteady flow modelling capabilities. 
Time integration is performed over a double periodic unit square domain $\Omega=\left(0,\,1\right)\times\left(0,\,1\right)$
by means of the BDF2 scheme. 
The initial horizontal and vertical velocity components are set as
\begin{equation*}
    u   =   \begin{cases}
                \tanh{\left(\dfrac{y-0.25}{\xi}\right)}    &   y\leq 0.5 \\
                \tanh{\left(\dfrac{0.75-y}{\xi}\right)}    &   y> 0.5
            \end{cases},
    \quad \text{and} \quad v   =   \delta \sin{\rbrackets{2\pi x}},
\end{equation*}
respectively. 
The pressure field is uniform and the free parameters are set as 
$\xi   =   1/30$ and $\delta  =   1/20$.

We consider $k=\left\lbrace 1,2,3,4\right\rbrace$ \hhohll formulations on a 
$h$-refined regular quadrilateral elements mesh sequence such that $\card{\Th} = 64 \ast 4^i, \, i=0,1,2,3,4$.
Time integration is carried out in the time interval $(t_0 = 0, t_F = 2]$ utilizing a constant time step $\Delta t = 10^{-3}$.
Upon time integration completion, the relative kinetic energy error, defined as 
\begin{equation*}
\dfrac{\mathcal{K}(t_{0})  -\mathcal{K}(t_{F})}{\mathcal{K}(t_{0})},
\end{equation*} 
is a measure of the amount of numerical dissipation introduced by the scheme.
Due to its lack of robustness in the inviscid limit, \hhohdiv fails to complete the test case,
regardless of the polynomial degree and the mesh density, the simulation blows up at early stages.
Accordingly only \hhohll results are presented hereafter.

In Figure~\ref{fig:double_shear_layer_convergence_cost} the relative kinetic energy error is plotted against mesh spacing and DOFs.
As expected, finer meshes and higher polynomial degrees provide increasingly small numerical dissipation.
Moreover, increasing the polynomial degree is always beneficial in terms of error versus DOFs.
Considering the finest two grids of the $h$-refined mesh sequence, 
the kinetic energy error convergence rates are $\left\lbrace 1.76,\, 2.13,\, 2.73,\, 2.71 \right\rbrace$ for $k=\left\lbrace 1,\,2,\,3,\,4 \right\rbrace$, respectively.
Probably, due to the high accuracy required to capture all the flow features, finer meshes and smaller time-steps would be required to reach asymptotic convergence rates. 
Indeed, as depicted in Figure~\ref{fig:double_shear_layer_vorticity} in case of $k=4$ \hhohll formulations, only the finest mesh allows to satisfactorily represent the tiniest vortical features.
This behaviour testifies about the challenges involved in the simulation of the double shear layer.
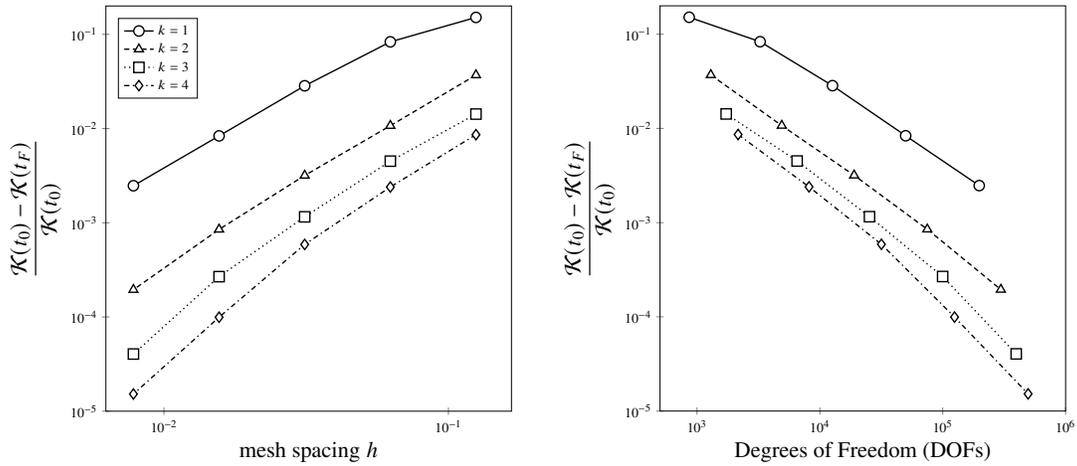
\begin{figure} [!t]
\centering
\begin{tikzpicture}[scale=0.65]
\begin{axis}    [
                unbounded coords=jump,
                ymode=log,xmode=log,
                scale only axis=true,
                width=0.5\textwidth,
                height=0.5\textwidth,
                xlabel=mesh spacing $h$,
                ylabel= $\dfrac{\mathcal{K}(t_{0})  -\mathcal{K}(t_{F})}{\mathcal{K}(t_{0})}$,
                label style={anchor=near ticklabel, font=\LARGE},
                ticks=major,
                tick pos=left,
                tick align=center,
                xmin=6.25E-3,
                xmax=1.666E-1,
                ymin=1E-5,
                ymax=2E-1,
                legend style={cells={anchor=west}},
                legend pos=north west,
                ]
\addplot    [
            thick,
            mark = *,
            mark size=3pt,
            mark options={black, fill=white}
            ]
            table[
            x expr = {\thisrowno{0}},
            y expr = {\thisrowno{2}},
            header=false]
            {DSL_k_p1.dat};
\addplot    [
            thick,
            densely dashed,
            mark = triangle*,
            mark size=3pt,
            mark options={black, solid, fill=white}
            ]
            table[
            x expr = {\thisrowno{0}},
            y expr = {\thisrowno{2}},
            header=false]
            {DSL_k_p2.dat};
\addplot    [
            thick,
            dotted,
            mark = square*,
            mark size=3pt,
            mark options={black, solid, fill=white}
            ]
            table[
            x expr = {\thisrowno{0}},
            y expr = {\thisrowno{2}},
            header=false]
            {DSL_k_p3.dat};
\addplot    [
            thick,
            dashdotted,
            mark = diamond*,
            mark size=3pt,
            mark options={black, solid, fill=white}
            ]
            table[
            x expr = {\thisrowno{0}},
            y expr = {\thisrowno{2}},
            header=false]
            {DSL_k_p4.dat};
\legend {$k=1$, $k=2$, $k=3$, $k=4$}
\end{axis}
\end{tikzpicture}
\hspace{5mm}
\begin{tikzpicture}[scale=0.65]
\begin{axis}    [
                unbounded coords=jump,
                ymode=log,xmode=log,
                scale only axis=true,
                width=0.5\textwidth,
                height=0.5\textwidth,
                xlabel=Degrees of Freedom (DOFs),
                ylabel= $\dfrac{\mathcal{K}(t_{0})  -\mathcal{K}(t_{F})}{\mathcal{K}(t_{0})}$,
                label style={anchor=near ticklabel, font=\LARGE},
                ticks=major,
                tick pos=left,
                tick align=center,
                xmin=5E+2,
                xmax=1E+6,
                ymin=1E-5,
                ymax=2E-1,
                legend style={cells={anchor=west}},
                legend pos=north east,
                ]
\addplot    [
            thick,
            mark = *,
            mark size=3pt,
            mark options={black, fill=white}
            ]
            table[
            x expr = {\thisrowno{1}},
            y expr = {\thisrowno{2}},
            header=false]
            {DSL_k_p1.dat};
\addplot    [
            thick,
            densely dashed,
            mark = triangle*,
            mark size=3pt,
            mark options={black, solid, fill=white}
            ]
            table[
            x expr = {\thisrowno{1}},
            y expr = {\thisrowno{2}},
            header=false]
            {DSL_k_p2.dat};
\addplot    [
            thick,
            dotted,
            mark = square*,
            mark size=3pt,
            mark options={black, solid, fill=white}
            ]
            table[
            x expr = {\thisrowno{1}},
            y expr = {\thisrowno{2}},
            header=false]
            {DSL_k_p3.dat};
\addplot    [
            thick,
            dashdotted,
            mark = diamond*,
            mark size=3pt,
            mark options={black, solid, fill=white}
            ]
            table[
            x expr = {\thisrowno{1}},
            y expr = {\thisrowno{2}},
            header=false]
            {DSL_k_p4.dat};
\end{axis}
\end{tikzpicture}
\caption{Double shear layer. \hhohll kinetic energy relative error at the final time. \emph{Left and right}: in dependence of the mesh spacing and of the DOFs of the discretized system, respectively.}\label{fig:double_shear_layer_convergence_cost}
\end{figure}
\begin{figure}[!t]
    \centering
	\includegraphics[width=0.4\textwidth]{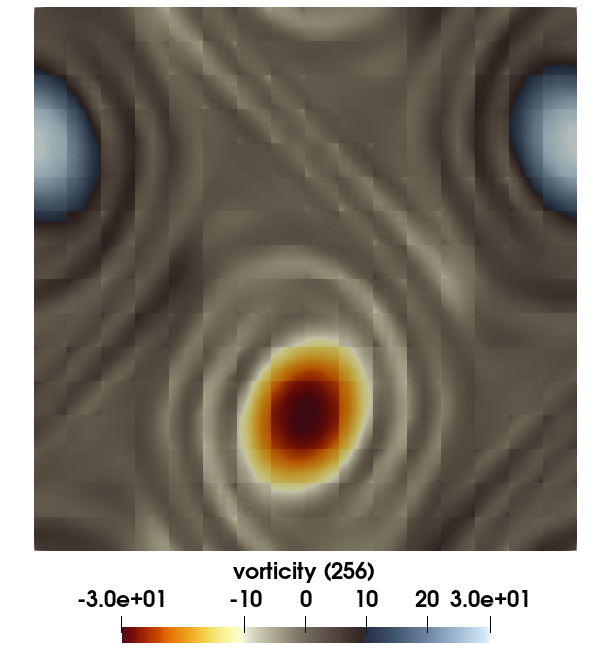}
	\includegraphics[width=0.4\textwidth]{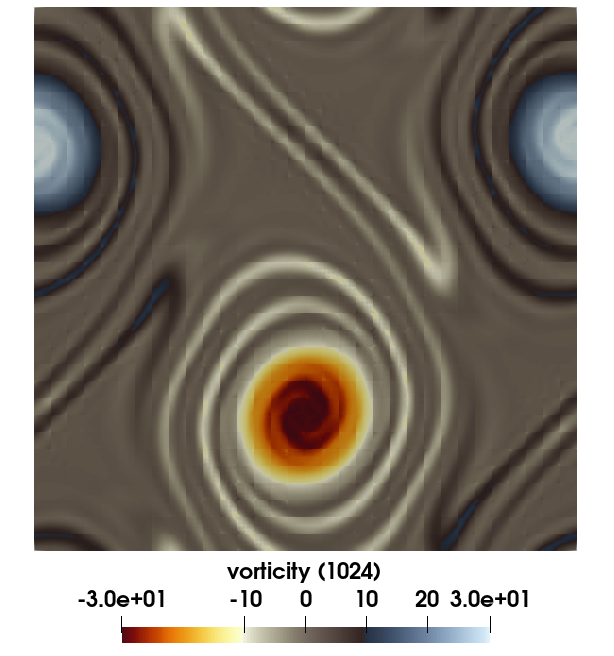}
	\includegraphics[width=0.4\textwidth]{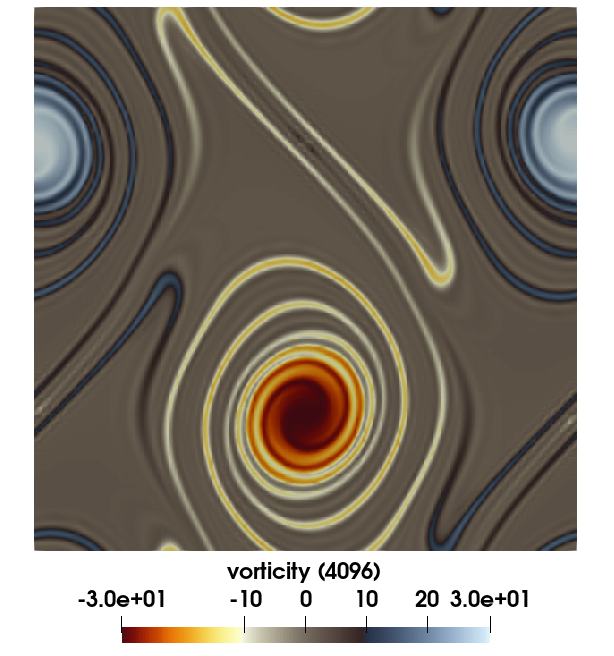}
	\includegraphics[width=0.4\textwidth]{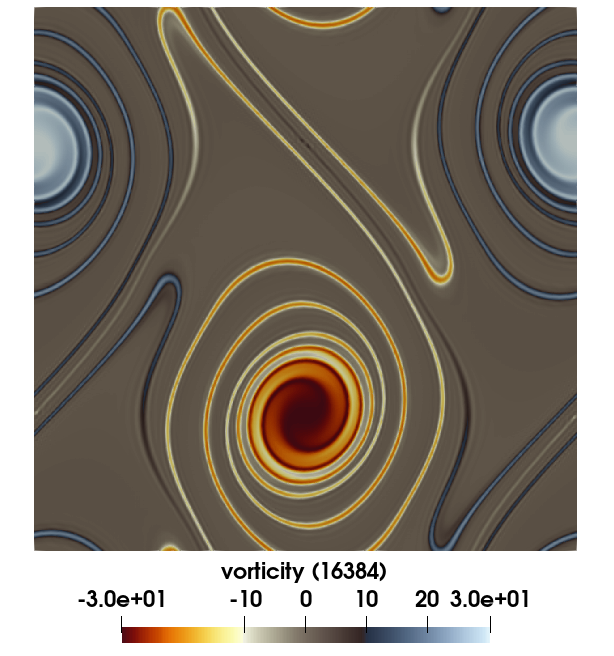}
    \caption{Double shear layer. Vorticity field at the final time obtained with $k=4$ \hhohll. The mesh cardinality, \ie 256, 1024, 4096 and 16384 elements, is reported within round brackets.}
	\label{fig:double_shear_layer_vorticity}
\end{figure}

\subsection{Lid-driven cavity flow}
\label{sec:LidDriven}
We consider the well-known lid-driven cavity flow problem over a unit square domain $\Omega=\left(0,\,1\right)\times\left(0,\,1\right)$ at Reynolds $\Reynolds=\frac{1}{\nu}=10^4$.
Dirichlet boundary conditions are imposed on the sliding top wall and the remaining stationary walls.
In order to evaluate the robustness of HHO formulations with respect to the mesh distortion
we consider a 12768 triangular elements mesh featuring anisotropic and stretched simplexes, see Figure~\ref{fig:lid_driven_cavity_mesh}.
$k=4$ \hhohll and \hhohdiv discretizations are applied for seeking steady state solutions starting from fluid at rest.
\begin{figure}[!t]
    \centering
    \includegraphics[width=.45\textwidth]{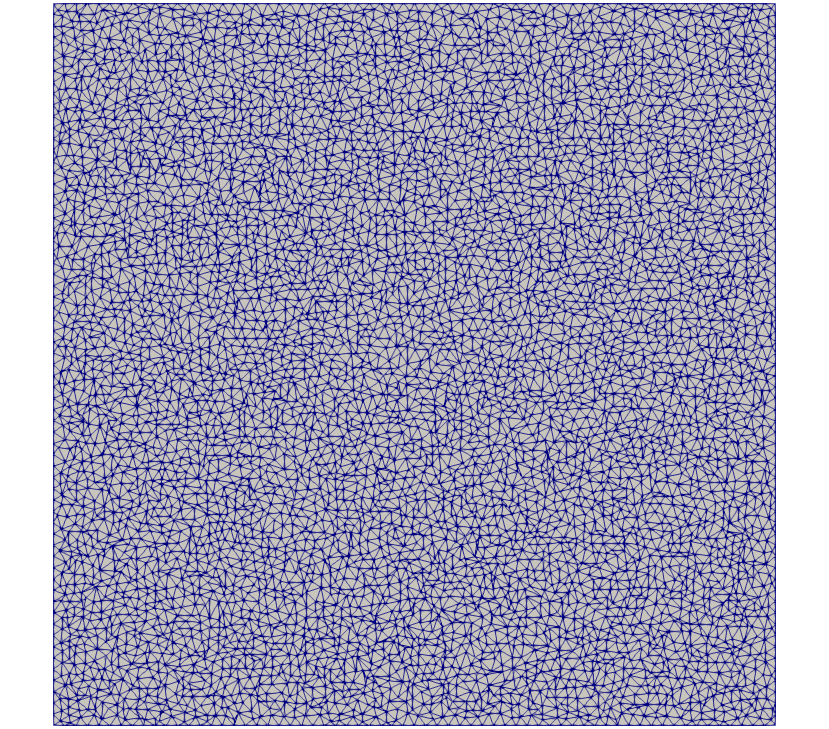}
    \hspace{5mm}
    \includegraphics[width=.433\textwidth]{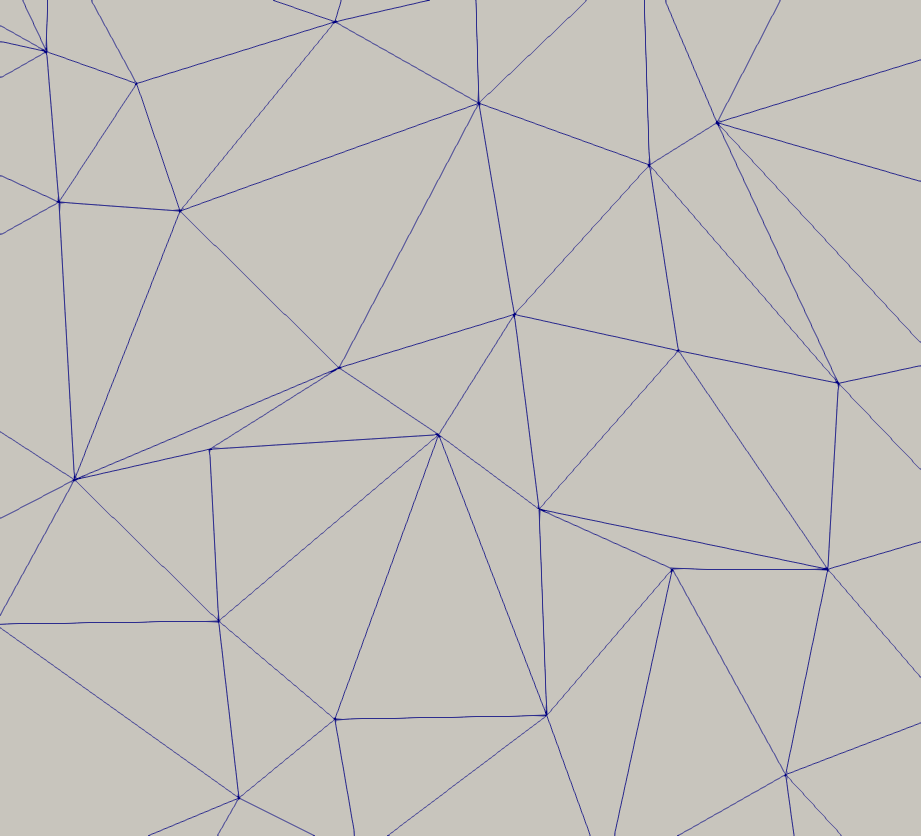}
    \caption{Lid-driven cavity. Computational mesh. \emph{Left and right}: global and detailed views, respectively.}
    \label{fig:lid_driven_cavity_mesh}
\end{figure}

In Figure~\ref{fig:lid_driven_cavity_velocity_profiles} the horizontal ($u$) and vertical ($v$) velocity components 
are plotted along vertical ($y=0.5$) and horizontal ($x=0.5$) centerlines of the cavity, respectively.
The two HHO formulation provide perfectly superimposed velocity profiles 
and are in very good agreement with reference solutions available from the literature \cite{Erturk.Corke.ea:05,Ghia.Ghia.ea:82}.
\begin{figure}\centering
\begin{tikzpicture}
  	[font=\footnotesize, spy using outlines={magnification=2.5, size=2.5cm, connect spies, fill=none}]
\begin{axis}   [height=9cm,
				width=9cm,
				xmin=-1,
				xmax=1,
				ymin=0,
				ymax=1,
        		xlabel={$u$},
        		ylabel={$y$},
        		legend style = { at={(0.205,0.975)}, anchor=north west, cells={anchor=west}},
        		axis x line*=top,
        		ytick pos=bottom]
\addplot[mark=o, only marks, thick, color=black!50]
			table [
			x expr = {\thisrowno{3}},
			y expr = {\thisrowno{0}},
			header=false]
			{LDC_erturk_y_vs_u.dat};
\addplot[mark=+, only marks, thick, color=black!75]
			table [
			x expr = {\thisrowno{3}},
			y expr = {\thisrowno{0}},
			header=false]
			{LDC_ghia_y_vs_u.dat};
\addplot   [only marks,
			mark=*,
            mark size=2pt,
            mark options={blue}]
			coordinates{(-1, -1)};  
\addplot   [only marks,
			mark=triangle*,
            mark size=2pt,
            mark options={orange}]
			coordinates{(-1, -1)};  
\addplot   [only marks,
			mark=*,
            mark size=0.5pt,
            mark options={blue}] 
            table[
            x expr = {\thisrowno{1}},
            y expr = {\thisrowno{0}},
            header=false] {LDC_hdiv_y_vs_u.dat};     
\addplot   [only marks,
			mark=triangle*,
            mark size=0.5pt,
            mark options={orange}] 
            table[
            x expr = {\thisrowno{1}},
            y expr = {\thisrowno{0}},
            header=false] {LDC_hll_y_vs_u.dat};
\legend{Erturk \ea \cite{Erturk.Corke.ea:05},Ghia \ea \cite{Ghia.Ghia.ea:82}, {$k=4$, \hhohdiv},{$k=4$, \hhohll}}
\end{axis}
\begin{axis}   [height=9cm, width=9cm, %
        		xmin=0,
        		xmax=1,
        		ymin=-1,
        		ymax=1,
        		xlabel={$x$},
        		ylabel={$v$},
        		axis y line*=right,
        		xtick pos=left]
\addplot[mark=o, only marks, thick, color=black!50]
			table [
			x expr = {\thisrowno{0}},
			y expr = {\thisrowno{3}},
			header=false]
			{LDC_erturk_x_vs_v.dat};
\addplot[mark=+, only marks, thick, color=black!75]
			table [
			x expr = {\thisrowno{0}},
			y expr = {\thisrowno{3}},
			header=false]
			{LDC_ghia_x_vs_v.dat};     
\addplot   [only marks,
			mark=*,
            mark size=0.5pt,
            mark options={blue}] 
            table[
            x expr = {\thisrowno{0}},
            y expr = {\thisrowno{1}},
            header=false] {LDC_hdiv_x_vs_v.dat};     
\addplot   [only marks,
			mark=triangle*,
            mark size=0.5pt,
            mark options={orange}] 
            table[
            x expr = {\thisrowno{0}},
            y expr = {\thisrowno{1}},
            header=false] {LDC_hll_x_vs_v.dat}; 
\end{axis}
\begin{scope}
    \spy on (6.85,1.85) in node [fill=none, anchor=south west] at (8.35,0.75);
    \spy on (5.50,6.85) in node [fill=none, anchor=south west] at (8.35,4.25);
    \spy on (0.55,5.25) in node [fill=none, anchor=south east] at (-0.75,4.25);
 	\spy on (2.25,0.55) in node [fill=none, anchor=south east] at (-0.75,0.75);
\end{scope}
\end{tikzpicture}
\caption{Lid-driven cavity at $\Reynolds = 10k$. Comparison with reference solutions based on horizontal and vertical velocity profiles over vertical and horizontal centerlines, respectively.}
\label{fig:lid_driven_cavity_velocity_profiles}
\end{figure}
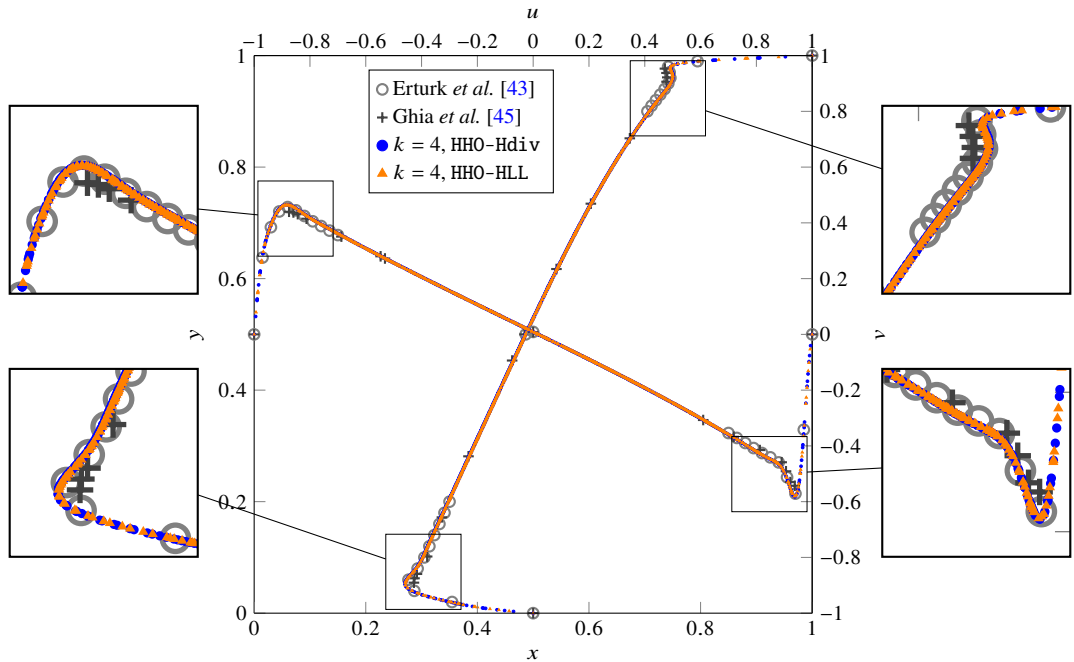

For the sake of comparison, in Figure~\ref{fig:lid_driven_cavity_warp velocity} 
we report the warp of the velocity magnitude fields obtained with \hhohdiv and \hhohll.
It is interesting to remark that the most noticeable difference between the two numerical solutions is related to the behavior at the top-left and top-right corners of the cavity,
where the velocity is discontinuous and the pressure gradient is steeper.
In particular, at the top-left corner, the \hhohdiv velocity solution shows a less good agreement with the weakly imposed no-slip boundary condition. 
\begin{figure}[!t]
    \centering
    \includegraphics[width=.45\textwidth]{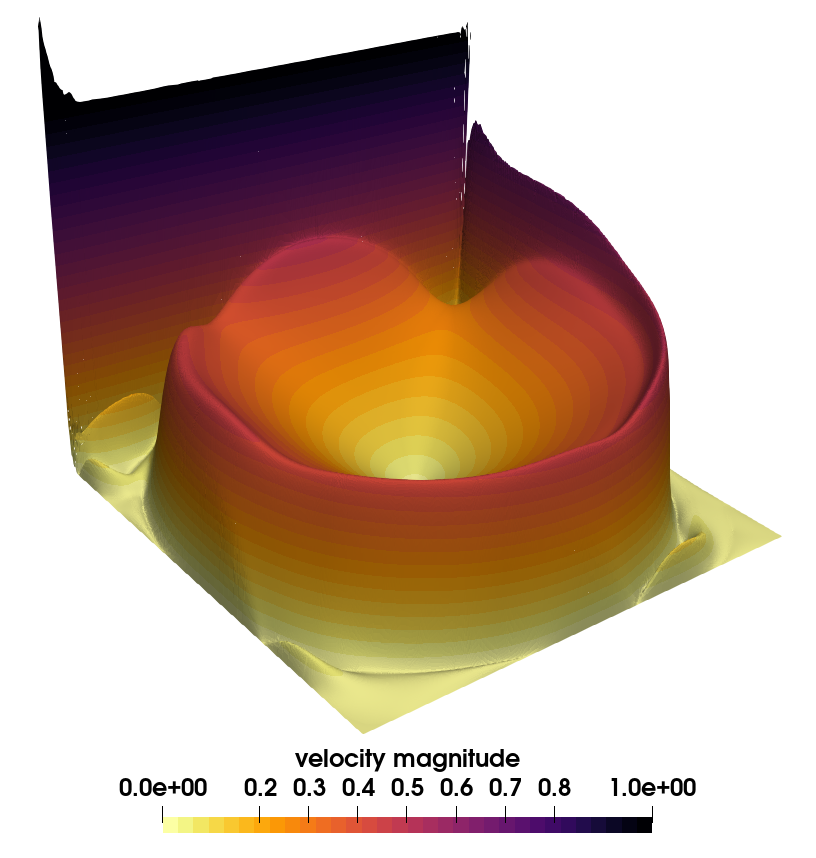}
    \hspace{5mm}
    \includegraphics[width=.45\textwidth]{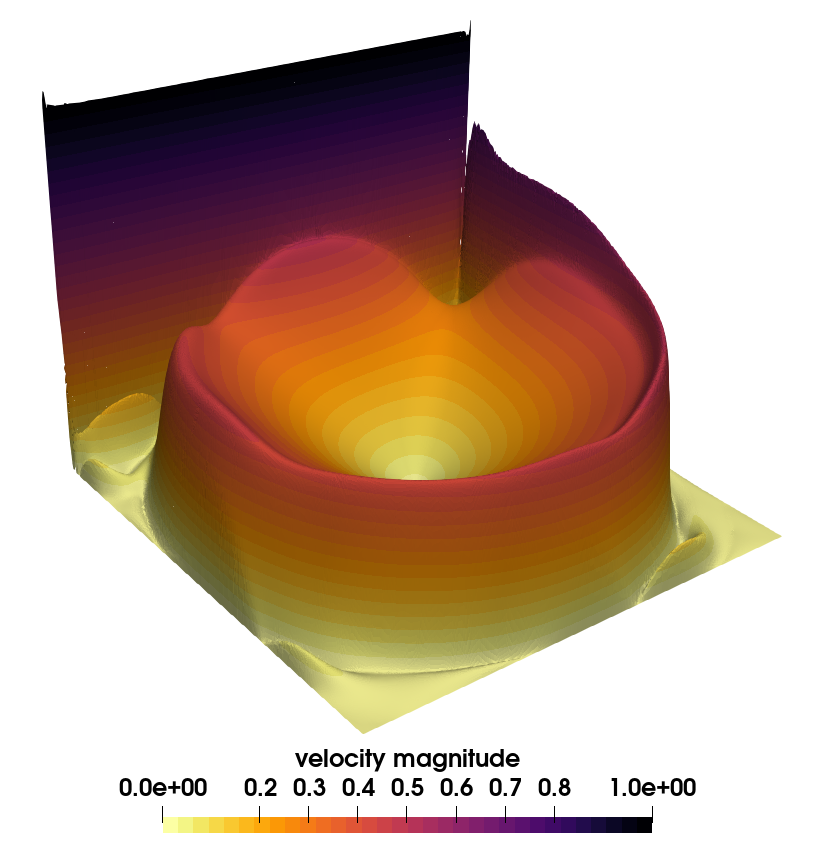}
    \caption{Lid-driven cavity at $\Reynolds = 10k$. Warp of the velocity magnitude field. \emph{Left and right}: \hhohdiv and \hhohll formulation, respectively.}
    \label{fig:lid_driven_cavity_warp velocity}
\end{figure}

\section{Conclusions}
We numerically validated two original Hybrid High-Order formulations designed for seeking approximate solutions of incompressible flow problems. 
The formulations allows to cope with Dirichlet and Neumann boundary conditions and, in particular, Dirichlet boundary conditions can be imposed weakly. 
From the accuracy viewpoint both formulations provides increased convergence rates with respect to discontinuous Galerkin discretizations in the diffusion dominated regime.
If we consider the polynomial degree leading the size of Jacobian matrix blocks as a reference, HHO $h$-convergence rates for both velocity and pressure unknowns are one order better than those of dG.
In case of \hhohdiv the improved convergence rates are maintained up to moderately high Reynolds numbers.

Besides providing numerical convergence rates, pressure-robustness and robustness in the inviscid limit are investigated performing especially conceived test cases. 
The \hhohll scheme, thanks to the introduction of Godunov fluxes based on HHL Riemann Solvers, can be employed for seeking approximate solution of the incompressible Euler equations. 
The \hhohdiv scheme is pressure-robust, kinetic energy preserving and yields mass conservation up to machine precision but is not stable in the inviscid limit. 
Robustness with respect to mesh distortion and grading is demonstrated for both formulations by solving the lid-driven cavity problem at high-Reynolds number over randomly distorted triangular elements meshes.

In conclusion \hhohdiv is best suited for low-to-moderate Reynolds number flows, thanks to improved convergence rates and exact mass conservation,
while \hhohll is able to cope with moderate-to-high Reynolds number flows, hence rivaling dG discretizations based on Godunov fluxes. 
While in this work we focused on 2D test cases, the formulations are well suited to be applied in 3D. 
Future works tackling 3D computations are planned but will require further efforts for improving the efficiency of the solution strategy. 
Indeed, while 2D computations can be effectively carried out relying on direct solvers, state-of-the-art preconditioned iterative solvers are required for tackling $h$-refined 3D mesh sequences.

\section*{Declarations}
\noindent
The authors declare that no funds, grants, or other support were received during the preparation of this manuscript.
The authors declare that they have no conflict of interest.
Data will be made available from the corresponding author on reasonable request.

\bibliographystyle{plain}      
\bibliography{hll}   

\end{document}